\documentclass[10pt]{amsart}
\usepackage[pagebackref]{hyperref}

\renewcommand{\d}{\mathrm{d}}%exterior derivative
\renewcommand{\L}{\Lambda}
\newcommand{\ts}{\textstyle }

\newcommand{\bbR}{{\mathbb R}}
\newcommand{\bbC}{{\mathbb C}}
\newcommand{\bbD}{{\mathbb D}}

\newcommand{\bbP}{{\mathbb P}}

\newcommand{\bbZ}{{\mathbb Z}}

\newcommand{\E}{{\mathrm e}}
\newcommand{\iC}{{\mathrm i}}

\newcommand{\SL}{\operatorname{SL}}
\newcommand{\PSL}{\operatorname{PSL}}
\newcommand{\SO}{\operatorname{SO}}

\newcommand{\SU}{\operatorname{SU}}

\newcommand{\Ric}{\operatorname{Ric}}

\newcommand{\Lie}{\operatorname{\mathsf{L}}}

\newcommand{\p}{\partial}

\newcommand{\cI}{\mathcal{I}}

\newcommand{\Cs}{\mathsf{C}}

\newcommand{\eug}{\operatorname{\mathfrak{g}}}

\newcommand{\eusu}{\operatorname{\mathfrak{su}}}

\newcommand{\bj}{{\bar{\jmath}}}

\newcommand{\w}{{\mathchoice{\,{\scriptstyle\wedge}\,}%
{{\scriptstyle\wedge}}{{\scriptscriptstyle\wedge}}%
{{\scriptscriptstyle\wedge}}}}
\newcommand{\lhk}{\mathbin{\hbox{\vrule height1.4pt width4pt depth-1pt 
             \vrule height4pt width0.4pt depth-1pt}}}

\newcommand{\be}{\begin{equation}}
\newcommand{\ee}{\end{equation}}
\newcommand{\bpm}{\begin{pmatrix}}
\newcommand{\epm}{\end{pmatrix}}

\numberwithin{equation}{section}

\newtheorem{proposition}{Proposition}

\theoremstyle{remark}
\newtheorem{definition}{Definition}

\newtheorem{remark}{Remark}
\newtheorem{example}{Example}

\begin{document}

\author[R. Bryant]{Robert L. Bryant}
\address{Duke University Mathematics Department\\
         P.O. Box 90320\\
         Durham, NC 27708-0320}
\email{\href{mailto:bryant@math.duke.edu}{bryant@math.duke.edu}}
\urladdr{\href{http://www.math.duke.edu/~bryant}%
         {http://www.math.duke.edu/\lower3pt\hbox{\symbol{'176}}bryant}}

\title[Hypersurfaces in unimodular complex surfaces]
      {Real hypersurfaces\\
       in\\
       unimodular complex surfaces}

\date{July 27, 2004}

\begin{abstract}
A unimodular complex surface is a complex $2$-manifold~$X$
endowed with a holomorphic volume form~$\Upsilon$.  
A strictly pseudoconvex real hypersurface in~$X$ inherits 
not only a CR-structure but a canonical coframing as well.

In this article, this canonical coframing is defined,
its invariants are discussed and interpreted
geometrically, and its basic properties are studied.  

A natural evolution equation for strictly
pseudoconvex real hypersurfaces in unimodular complex surfaces
is defined, some of its properties are discussed,
and several examples are computed.

It is shown that a real-analytic $3$-manifold endowed with
a real-analytic coframing satisfying the structure equations 
can be real-analytically embedded as a pseudoconvex hypersurface 
in a unimodular complex surface in such a way that the 
induced canonical coframing is the given one.  Moreover,
this embedding is essentially unique up to unimodular
biholomorphism.

The locally homogeneous examples are determined 
and used to illustrate various features of the geometry
of the induced structure on the hypersurface.

The invariants of the underlying CR-structure 
are expressed in terms of the invariants of the coframing.
\end{abstract}



\thanks{
Thanks to 
Duke University for its support via a research grant, 
to the NSF for its support via DMS-0103884, 
and to Columbia University for its support 
via an Eilenberg Visiting Professorship.  
\hfill\break
\hspace*{\parindent} 
This is Version~$1.0$ of CR3mfds.tex
}

\maketitle

\section{Introduction}\label{sec: intro}

The purpose of this article is to define and study 
a flow for nondegenerate real hypersurfaces in a complex 
surface endowed with a holomorphic volume form.  This 
flow is invariant under unimodular biholomorphisms
and is weakly parabolic in an appropriate sense.  
Neither uniqueness nor short time existence for this 
flow has yet been proved and it should be an interesting
problem to do so.  

The remainder of this introduction will be a guide to
the sections of the article.

A unimodular complex surface is a complex $2$-manifold~$X$
endowed with a holomorphic volume form~$\Upsilon$.  
A strictly pseudoconvex real hypersurface in~$X$ inherits 
not only a CR-structure but a canonical coframing as well.
This canonical coframing is defined in Section~\ref{sec: Rhypersurfaces}
by its structure equations (cf. Proposition~\ref{prop: cancoframing}).

It is worth emphasizing that the ambient holomorphic volume
form is needed in order to define this coframing: There is no 
canonical coframing of strictly pseudoconvex real hypersurfaces 
in~$\bbC^2$ that is invariant under the full pseudogroup of 
biholomorphisms of~$\bbC^2$.

Specifically, it is shown that if~$M\subset X$ is 
a strictly pseudoconvex real hypersurface in~$X$, 
then the pullback of~$\Upsilon$ to~$M$ can be written
uniquely in the form~$M^*\Upsilon=\theta\w\eta$ for
$1$-forms~$\theta$ and~$\eta$ on~$M$ satisfying
\be
\theta = \bar\theta
\qquad\text{and}\qquad
\d\theta = \iC\,\eta\w\overline{\eta}\,.
\ee
The~$1$-form~$\theta$ depends on two derivatives of
a local defining function for~$M$ while the $1$-form~$\eta$
depends on three derivatives of such a local defining function.

It is then shown that there exist functions~$a = \bar a$ 
and~$b$ on~$M$ such that
\be
\d\eta = 2\iC\,\theta\w\bigl(a\,\eta + b\,\overline{\eta}\,\bigr).
\ee
The functions~$a$ and~$b$ are the \emph{primary invariants}
of~$M$.  They are the fourth order invariants of~$M$ as a
hypersurface in~$(X,\Upsilon)$: Given any point~$p\in M$, 
there exist $p$-centered holomorphic coordinates~$(w,z)$ on a
$p$-neighborhood~$U\subset X$ such that~$U^*\Upsilon = \d w\w\d z$
and so that~$M\cap U$ is defined in~$U$ by an equation of
the form
\be
\mathrm{Im}(w) = 
{\ts\frac12}\,z\bar z\,\bigl(1 + b(p)\,z^2 
                 + {\ts\frac32} a(p)\,z\bar z 
                 +  \overline{b(p)}\,{\bar z}^2 \bigr)
    + R_5\bigl(z,\mathrm{Re}(w)\bigr)
\ee
where~$R_5$ is a function on a neighborhood of~$(0,0)\in\bbC\times\bbR$
that vanishes to order~$5$ at~$(z,u)= (0,0)$.

The invariants $a$ and~$b$ are then discussed 
and interpreted geometrically, and their 
basic properties are studied.  For example, 
$a$ plays the role of `mean curvature' 
in the sense that its vanishing is the 
Euler-Lagrange equation for the functional
defined by integrating the canonical volume 
form~$\frac12\,\theta\w\d\theta$ over~$M$.
The condition for CR-flatness is expressed in
terms of the invariants~$a$ and~$b$ and their
derivatives.  It is no surprise that there should 
be such an expression since~$a$ and~$b$, together
with the canonical coframing~$(\theta,\eta)$, is
a \emph{complete} set of invariants (in Cartan's sense) 
for nondegenerate real hypersurfaces in unimodular surfaces.%
\footnote{Roughly speaking, this means that any 
pointwise differential invariant for such hypersurfaces
can be expressed in terms of~$a$ and~$b$ and their
derivatives with respect to the coframing~$(\theta,\eta)$.}

A nondegenerate hypersurface~$M\subset X$
is locally homogeneous under the unimodular
biholomorphism pseudogroup if and only if~$a$ and~$b$
are constant functions.   It is shown that
there is a $3$-parameter family of examples
of locally homogeneous hypersurfaces and these
examples are constructed explicitly. 

In Section~\ref{ssec: embeddingresults},
it is shown that a real-analytic $3$-manifold endowed with
a real-analytic coframing satisfying the structure equations
of Proposition~\ref{prop: cancoframing}
can be real-analytically embedded as a pseudoconvex hypersurface 
in a unimodular complex surface in such a way that the 
induced canonical coframing is the given one.  Moreover,
this embedding is essentially unique up to unimodular
biholomorphism.  This shows that the coframing on the
hypersurface captures all of the ambient local unimodular
biholomorphism invariants. 

In Section~\ref{ssec: flows}
a natural evolution equation for strictly
pseudoconvex real hypersurfaces in unimodular complex surfaces
is defined, some of its properties are discussed,
and some homogeneous examples are computed.  The basic
idea is this:  If~$T$ is the Reeb vector field on~$M$
associated to the contact form~$\theta$, then~$\iC T$
is a canonical normal vector field along~$M\subset X$.
The \emph{unimodular normal flow} flow is then defined 
so that~$\iC T$ is the velocity field of the flow.

This flow is somewhat subtle.  The normal vector 
field~$\iC T$ depends on third order information at
each point of the hypersurface~$M$, and the flow
preserves the contact structure on~$M$ (i.e., $\theta$
changes only by a multiple).  However, if one is
allowed to reparametrize~$M$ during the flow (in
a manner that is not invariant under the unimodular
biholomorphism pseudogroup), then one can replace the
flow by a second order flow, one whose linearization
is subelliptic and that thus, presumably, 
has short-time existence and uniqueness 
for smooth initial data.

Explicitly, when~$M$ is expressed in local unimodular
coordinates~$(w,z)$ as a graph of the form~$\mathrm{Im}(w)
= F_0\bigl(z,\mathrm{Re}(w)\bigr)$, the evolution equation
for the function~$F=F(t,z,u)$ that is induced by the 
unimodular normal flow is
\be
F_t = \Bigl[2\left[
         (1{+}{F_u}^2)F_{z\bar z}
         - \iC (1{-}\iC\,F_u)F_{\bar z}F_{uz}
         + \iC (1{+}\iC\,F_u)F_{z}F_{u\bar z}
          + F_zF_{\bar z}F_{uu} \right]\Bigr]^{1/3},
\ee
with initial condition~$F(0,z,u) = F_0(z,u)$. 
(The quantity inside the radical on the right hand side 
of this equation is always nonzero for a pseudoconvex
hypersurface represented as a graph in this way.)

As a first step in understanding this flow, 
the evolution equations for the invariants are computed
and examples of the flow are computed for homogeneous hypersurfaces.
In most cases, the flow exists only for 
a finite time (either in the future or the past) 
and the nature of the singularity that develops is discussed.
An example shows that this flow does not preserve CR-flatness,
in general.

Section~\ref{sec: CRinvs} is provided for the reader's 
convenience.  It is an account (closely following
Cartan's original article) of the equivalence method
development of the Cartan connection associated to a 
pseudoconvex real hypersurface in a complex surface.  
This is used both as a contrast to the unimodular case 
and as a reference for the notions from CR-geometry 
that are used in the rest of the article.  

Section~\ref{sec: generalizations} outlines the changes
needed to generalize these constructions to the case of
a pseudoconvex hypersurface~$M^{2n+1}$ in a higher dimensional 
unimodular complex manifold of dimension~$n{+}1$ 
and to the case in which the ambient complex manifold 
is endowed with a real volume form (rather than a 
holomorphic volume form).  It turns out, in either case, 
that one can still define a canonical~$\SU(n)$-structure 
on~$M$ and a flow.

Finally, it is a pleasure to thank Richard Hamilton, 
whose questions about a different flow for CR-hypersurfaces 
in K\"ahler manifolds that was studied by Huisken and 
Klingenberg~\cite{MR01f:53141,MR03d:32046} led me to search
for a flow invariant under the unimodular biholomorphism group.

\setcounter{tocdepth}{2}
\tableofcontents

\section{Real hypersurfaces in unimodular surfaces}
\label{sec: Rhypersurfaces}

Let~$X$ be a complex $2$-manifold and let~$\Upsilon$ 
be a nowhere-vanishing holomorphic $2$-form on~$X$.
The pair~$(X,\Upsilon)$ will be said to be a \emph{unimodular surface}.

A coordinate chart~$z = (z^1,z^2):U\to \bbC^2$ 
on an open set~$U\subset X$ will be said to be 
\emph{unimodular} if the pullback of~$\Upsilon$
to~$U$ is equal to~$\d z^1\w \d z^2$.  

A unimodular surface  has an atlas of unimodular coordinate charts.
In fact, if~$z^1:U\to\bbC$ is a holomorphic function on~$U\subset X$
and~$\d z^1$ is nonvanishing at~$p\in U$, then there is a 
$p$-neighborhood~$V\subset U$ on which there exists another
holomorphic function~$z^2:V\to\bbC$ such that~$z = (z^1,z^2):V\to\bbC^2$
is a unimodular coordinate chart.  Note that, if the fibers of~$z^1:V\to\bbC$
are connected, then the function~$z^2$ is determined up to the addition 
of a holomorphic function of~$z^1$.

\subsection{Induced geometry on hypersurfaces}
\label{ssec: inducedgeometry}
Now let~$M\subset X$ be a real hypersurface.  
For simplicity,~$M$ will be assumed to be smooth, 
although weaker differentiability hypotheses 
would suffice for most of the constructions.

Let~$\Phi$ be the pullback of~$\Upsilon$ to~$M$.
Since~$\SL(2,\bbC)$ acts transitively on the set 
of real hyperplanes in~$\bbC^2$, 
the algebraic type of~$\Phi$ does not vary on~$M$.  

In particular, any point~$p\in M$
will have an open neighborhood~$V\subset M$ on which
there will exist $1$-forms~$\theta = \overline\theta$
and~$\eta$ such that~$V^*\Phi = \theta\w\eta$ 
and where~$\theta\w\eta\w\overline{\eta}\not=0$.
These forms are not unique, but any other pair~$(\theta^*,\eta^*)$
on~$V$ with these properties will satisfy
\be\label{eq: 0coframetransition}
\bpm \theta^*\\ \eta^*\epm
= \bpm \lambda & 0 \\ \mu & \lambda^{-1}\epm \bpm \theta\\ \eta \epm
\ee
for some functions~$\lambda=\bar \lambda\not=0$ and $\mu$ on~$V$.

In particular, the equation~$\theta=0$ defines 
a global $2$-plane field~$D\subset TM$ 
that is independent of the local choice of~$\theta$.%
\footnote{The $2$-plane field~$D$ is known as the 
\emph{complex tangent space} of~$M$ in the CR literature.}
Since~$\theta^*\w\d\theta^* = \lambda^2\,\theta\w\d\theta$,
whether or not~$\theta\w\d\theta$ vanishes on~$V$ 
is independent of the choice of coframing~$(\theta,\eta)$.

\begin{definition}
A real hypersurface~$M\subset X$ is said to be \emph{nondegenerate}
if its $2$-plane field~$D\subset TM$ is contact, i.e., 
if $\theta\w\d\theta$ is nowhere-vanishing for any locally 
defined $1$-forms~$\theta=\overline{\theta}$ and~$\eta$
satisfying~$\Phi = \theta\w\eta$.
\end{definition}

\begin{remark}[Pseudoconvexity]
The condition that is being called `nondegenerate' here 
is also known as the condition of `strict pseudoconvexity' 
in the CR-literature.  This coincidence is special to hypersurfaces
in complex $2$-manifolds.  In higher dimensions, `nondegenerate'
is implied by `strict pseudoconvexity', but not \emph{vice versa}.
\end{remark}

\begin{proposition}[The canonical coframing]\label{prop: cancoframing}
Let~$M\subset X$ be a nondegenerate real hypersurface.  Then
there exist unique $1$-forms $\theta=\overline{\theta}$ 
and~$\eta$ on~$M$ such that
\begin{enumerate}
\item $\Phi = \theta \w \eta$, and
\item $\d\theta = \iC\,\eta\w \overline{\eta}$.
\end{enumerate}
\end{proposition}

\begin{proof}
Let~$V\subset M$ be an open set on which there 
exist $1$-forms~$\theta=\overline{\theta}$ and~$\eta$ 
satisfying~$V^*\Phi = \theta\w\eta$.  Then there exist
functions~$L=\overline{L}$ and~$P$ on~$V$ such that
\be
\d\theta = \iC L\,\eta\w\overline{\eta} 
+ \iC\,\bigl(P\,\eta - \bar P\,\overline{\eta}\bigr)\w\theta.
\ee
Since~$M$ is nondegenerate,~$\theta\w\d\theta$ is 
nonvanishing on~$V$.  Thus,~$L$ is nowhere-vanishing on~$V$.

Now let~$\theta^* = \lambda\,\theta$ and~$\eta^* = \mu\,\theta 
+ \lambda^{-1}\,\eta$ 
for some functions~$\lambda=\bar\lambda\not=0$ and~$\mu$ on~$V$.
Then there exist functions~$L^* = \overline{L^*}\not=0$ and~$P^*$ on~$V$
such that
\be
\d\theta^* = \iC L^*\,\eta^*\w\overline{\eta^*} 
+ \iC\,\bigl(P^*\,\eta^* - \overline{P^*}\,\overline{\eta^*}\bigr)\w\theta^*.
\ee

Since~$\theta^*$ is a multiple of~$\theta$, computation yields
\be
\left.
\begin{aligned}
\d\theta^* 
&\equiv \lambda\,\d\theta 
 \equiv \lambda\,\iC L\,\eta\w\overline{\eta}\\
&\equiv  \lambda\,\iC L\,(\lambda\eta^*)\w\overline{(\lambda\eta^*)}
 \equiv  \iC\,(\lambda^3L)\,\eta^*\w\overline{\eta^*}\quad
\end{aligned}\right\}\ \mod \theta\ (\equiv\mod\theta^*).
\ee
Thus,~$L^* = \lambda^3\, L$.  In particular,  
there is a unique choice of~$\lambda$, 
namely~$\lambda= L^{-1/3}$, such that~$L^*=1$. 
 
It follows that the equation
\be\label{eq: 1streqL=1}
\d\theta = \iC\,\eta\w\overline{\eta} 
+ \iC\,\bigl(P\,\eta - \bar P\,\overline{\eta}\bigr)\w\theta
\ee
uniquely defines~$\theta$ on~$V$ and therefore globally defines~$\theta$
on~$M$.  Consequently,~$\eta$ is defined on~$M$ up to the addition
of a multiple of~$\theta$.

However, replacing~$\eta$ by~$\eta + \bar P\,\theta$ 
in the coframing shows that one can arrange
\be\label{eq: etafixed}
\d\theta = \iC\,\eta\w\overline{\eta}.
\ee
Moreover, once~$\theta$ is fixed, 
~\eqref{eq: etafixed} uniquely determines~$\eta$ globally on~$M$.
\end{proof}

\begin{remark}[The torsion tensor]\label{rem: torten}
As the construction in the proof demonstrates, two coframings
related as in~\eqref{eq: 0coframetransition}, satisfy
\be
L^*\bigl(\iC\,\theta^*\w\eta^*\w\overline{\eta^*}\bigr)^3
= L\,\bigl(\iC\,\theta\w\eta\w\overline{\eta}\bigr)^3,
\ee
so the local expression
\be
L\bigl(\iC\,\theta\w\eta\w\overline{\eta}\bigr)^3
= (\theta\w\d\theta)\circ\bigl(\iC\,\theta\w\eta\w\overline{\eta}\bigr)^2
\ee
is the restriction to the domain of the coframing of a well-defined
global section~$\Cs$ of the line bundle~$\bigl(\L^3(T^*M)\bigr)^{\otimes 3}$.
The nondegeneracy condition is just that this section~$\Cs$ 
be nowhere-vanishing.
\end{remark}

\begin{remark}[Coordinate expressions]\label{rem: coordexprs}
It is useful to look at the tensor~$\Cs$ in local unimodular coordinates.
Suppose that~$p$ is a point of~$M$ and let~$(w,z):U\to\bbC^2$ be
a unimodular coordinate system centered on~$p$ such that the
hypersurface~$\textrm{Im}(w)=0$ is tangent to~$M$ at~$p$.  
Writing~$w = u + \iC v$, the hypersurface~$M$ can be described
near~$p$ in these coordinates by an equation of the form
\be\label{eq: definingF}
v = F(u,z,\bar z)
\ee
where~$F$ is a real-valued function that vanishes at~$(0,0,0)$
and has its first derivatives vanishing there.

One can then write~$\Phi = M^*(\Upsilon) = \theta_0\w\eta_0$ where
\be
\theta_0 = \bar\theta_0
       = \d u - \frac{\iC\,F_z}{1-\iC\,F_u}\,\d z
              + \frac{\iC\,F_{\bar z}}{1+\iC\,F_u}\,\d{\bar z}
\qquad\text{and}\qquad
\eta_0 = (1{+}\iC\,F_u)\,\d z.
\ee
Thus,~$\iC\,\theta_0\w\eta_0\w\overline{\eta_0} 
= \iC\,(1+{F_u}^2)\,\d u \w\d z\w\d{\bar z}$.
Computation yields
\be
\Cs = L_0\,\bigl(\iC\,\theta_0\w\eta_0\w\overline{\eta_0}\bigr)^3 
    = \bigl(1 + {F_u}^2\bigr)^3L_0\,
            \bigl(\iC\,\d u \w\d z\w\d{\bar z}\bigr)^3,
\ee
where
\be
L_0
= 2\,\frac{(1{+}{F_u}^2)F_{z\bar z} + F_zF_{\bar z}F_{uu}
   - \iC (1{-}\iC\,F_u)F_{\bar z}F_{uz}
   + \iC (1{+}\iC\,F_u)F_{z}F_{u\bar z}}{\bigl(1 + {F_u}^2\bigr)^3}.
\ee

Thus, nondegeneracy at~$p$ is equivalent to~$F_{z\bar z}(0,0,0)\not=0$
(since the first derivatives of~$F$ vanish at~$(0,0,0)$).

This computation shows that the normalized~$\theta 
= (L_0)^{-1/3}\,\theta_0$ of Proposition~\ref{prop: cancoframing} 
depends on two derivatives of the defining function~$F$.  
Moreover, the normalized~$\eta$ depends on three derivatives of~$F$.
\end{remark}

\subsubsection{Invariants in the nondegenerate case}
Let~$M\subset X$ be a nondegenerate hypersurface and let~$(\theta,\eta)$
be the canonical coframing described in Proposition~\ref{prop: cancoframing}.

\textsl{Primary invariants.}
Since~$\d(\theta\w\eta) = \theta\w\d\eta = 0$, it follows that there
exist functions~$a$ and~$b$ on~$M$ such that
\be\label{eq: deta}
\d\eta = 2\iC\,\theta\w(\,a\,\eta + b\,\overline{\eta}\,).
\ee
(The introduction of the coefficient~$2\iC$ in this equation 
simplifies later formulae.) Since
\be
0 = \d(\d\theta) = \d(\iC\,\eta\w\overline{\eta})
= 2(\bar a - a)\,\theta\w\eta\w\overline{\eta},
\ee
it follows that~$a = \bar a$ is real-valued.  
 
The functions~$a$ and~$b$ will be referred to as the 
\emph{primary invariants} of~$M$.  Since the coframing~$(\theta,\eta)$
depends on three derivatives of a defining function~$F$ 
as in~\eqref{eq: definingF}, 
it follows that~$a$ and~$b$ are fourth-order expressions in~$F$.

\begin{remark}[Scaling effects]\label{rem: scaleffects}
For any constant~$\lambda\in\bbC^*$, the pair~$(X,\lambda\Upsilon)$
defines a unimodular surface and this constant scaling will 
affect the canonical coframing.  

Explicitly, replacing~$\Upsilon$ by~$\lambda\,\Upsilon$ replaces
the canonical coframing~$(\theta,\eta)$ on a nondegenerate 
hypersurface~$M\subset X$ by the coframing
\be
\bigl(\,|\lambda|^{2/3}\theta,\,\lambda|\lambda|^{-2/3}\eta\,\bigr).
\ee
Hence, under this change, the primary invariants~$(a,b)$ are replaced by
\be
\bigl(\,|\lambda|^{-2/3}\,a,\,\lambda^{-2}|\lambda|^{4/3}\,b\,\bigr).
\ee

Finally, if one reverses the complex structure on~$X$,
effectively by replacing~$\Upsilon$ by~$\overline{\Upsilon}$, 
one finds that the canonical coframing~$(\theta,\eta)$ is
then replaced by~$(-\theta,-\overline{\eta})$, so that the
primary invariants~$(a,b)$ are replaced by~$(a,\overline{b})$.
\end{remark}

\textsl{Secondary invariants.}
Since~$\theta$, $\eta$, and~$\overline{\eta}$ are linearly
independent and so form a basis for the $1$-forms on~$M$,
there exist unique functions~$a_\theta= \overline{a_\theta}$, 
$a_\eta$, and $a_{\overline{\eta}}= \overline{a_\eta}$ 
satisfying
\be\label{eq: daexp}
\d a = a_\theta\,\theta + a_\eta\,\eta 
           + a_{\overline{\eta}}\,\overline{\eta}\,.
\ee
Similarly, there exist unique functions~$b_\theta$, $b_\eta$, 
and~$b_{\overline{\eta}}$ satisfying
\be\label{eq: dbexp}
\d b = b_\theta\,\theta + b_\eta\,\eta 
           + b_{\overline{\eta}}\,\overline{\eta}\,.
\ee
Taking the exterior derivative of~\eqref{eq: deta} yields
\be\label{eq: beta=aetabar}
b_\eta =  a_{\overline{\eta}}\ \ (\ = \overline{a_\eta}).
\ee
The quantities~$a_\theta$, $a_\eta$, $b_\theta$, 
and~$b_{\overline{\eta}}$ are the \emph{secondary invariants} 
of~$M$.

\textsl{Derivative commutants.}
In general, for a function~$f$ on~$M$, set
\be\label{eq: dfexpanded}
\d f = f_\theta\,\theta + f_\eta\,\eta + f_{\bar\eta}\,\overline{\eta}.
\ee
Of course, these differentiation operations are not independent.
For example, taking the exterior derivative of~\eqref{eq: dfexpanded}
yields the relations
\be\label{eq: cofderexchformulae}
\begin{aligned}
f_{\eta\bar\eta} - f_{\bar\eta\eta} &= \iC f_\theta\,,\\
f_{\theta\eta} - f_{\eta\theta} 
    &= \phantom{-}2\iC a\,f_{\eta} - 2\iC\overline{b}\,f_{\bar\eta}\,,\\
f_{\theta\bar\eta} - f_{{\bar\eta}\theta} 
    &= -2\iC a\,f_{\bar\eta} + 2\iC b\,f_\eta\,.
\end{aligned}
\ee

\textsl{A subelliptic Laplacian.}
The operator~$L:C^\infty(M)\to C^\infty(M)$ defined by
\be
L f = f_{\eta\bar\eta} + f_{\bar\eta\eta}
\ee
is a subelliptic real operator that plays the role of the
Laplacian in this geometry.

\subsubsection{Some interpretations of the invariants}
\label{sssec: Invinterps}
Let~$(\theta,\eta)$ be the canonical coframing of a nondegenerate
hypersurface~$M\subset X$ where~$(X,\Upsilon)$ is a unimodular 
complex surface.  The invariants found thus far have natural 
interpretations and lend themselves to various constructions.

For example, there is a natural volume form on~$M$, 
namely~$\theta\w\d\theta = \iC\theta\w\eta\w\overline{\eta}$. 

There is also a real vector field~$T$ on~$M$ defined by
the conditions~$\eta(T)=0$ and~$\theta(T)=1$.  The vector field~$T$
is the \emph{Reeb vector field} of the contact form~$\theta$.
The vector field~$\iC T$ is a canonical normal vector field along~$M$. 
It points towards the `pseudoconcave' side of the hypersurface, 
i.e., the side that contains the holomorphic discs whose boundaries 
lie in~$M$.

This intrinsic normal can be used to describe variations of 
nondegenerate hypersurfaces:  If~$\iota:(-\epsilon,\epsilon)\times M\to X$ 
is a compactly supported $1$-parameter family of immersions of~$M$ 
into~$X$ as nondegenerate hypersurfaces, then one can reparametrize
this family so that
\be
\iota_*\left(\frac{\partial\hfil}{\partial t}(t,p)\right)
= f(t,p)\,J\bigl(T(t,p)\bigr),
\ee
for~$(t,p)\in (-\epsilon,\epsilon)\times M$, i.e., 
so that the variation is `normal'.  
The function~$f:(-\epsilon,\epsilon)\times M\to\bbR$ 
then can be regarded as the `normal displacement function'.

In this case, one has
\be
\iota^*(\Upsilon) = (\theta + \iC f\,\d t) \w \eta
\ee
where~$(\theta,\eta)$ is the ($t$-dependent) induced canonical coframing
on~$M$.  Computation yields
\be
\d\theta 
= \iC\,\eta\w\overline{\eta} 
          +\bigl( {\ts\frac13}(4af+f_{\eta\bar\eta}+f_{\bar\eta\eta})\,\theta 
         + \iC(\,f_\eta\,\eta-f_{\bar\eta}\,\overline{\eta}\,)\bigr) \w \d t.
\ee

By this and integration by parts, 
if the variation is supported in a compact domain~$K\subset M$, 
the volume function
\be
V_K(t) = \int_{\{t\}\times K} \iC\,\theta\w\eta\w\overline{\eta}
\ee
satisfies
\be
V_K'(t) =  -\frac83\,\int_{\{t\}\times K} 
              fa\ \iC\,\theta\w\eta\w\overline{\eta}.
\ee
Since, for a given nondegenerate immersion~$\iota_0:M\to X$,
one can construct a family~$\iota:\bbR\times M\to X$
with any given compactly supported function~$f_0 = f(0,\cdot):M\to\bbR$,
one sees that $\iota_0:M\to X$ is volume critical for all compactly 
supported variations if and only if it satisfies $a_0\equiv0$.  
Thus, the invariant~$a$ can be regarded as an analog of `mean curvature'
of the nondegenerate immersion.  (Keep in mind
that~$a$ is a fourth order invariant of the immersion.  
This is not surprising, though, since the volume functional 
is second order.)

One interpretation of the invariant~$b$
has to do with the geometry of the surface~$\Sigma$ that is
the quotient of~$M$ by the flow of the Reeb vector field~$T$. 
Since~$\iC\,\eta\w\overline{\eta} = \d \theta$ is closed and
semi-basic for the projection~$M\to\Sigma$, it follows that it
is the pullback to~$M$ of a well-defined, nonvanishing $2$-form 
on~$\Sigma$.  However, the quadratic form~$\eta{\circ}\overline{\eta}$
satisfies
\be
\Lie_T\,\eta{\circ}\overline{\eta} 
= 2\iC\,\bigl(b\,{\overline{\eta}}^2- \overline{b}\,\eta^2\bigr).
\ee
Thus, the quadratic form~$\eta{\circ}\overline{\eta}$ is not the
pullback to~$M$ of a metric on~$\Sigma$ unless~$b\equiv0$.  

However, the equation~$b\equiv0$ implies, via~\eqref{eq: beta=aetabar},
that~$a_{\bar\eta} = a_\eta = 0$.  Then the exterior derivative
of the equation~$\d a = a_\theta\,\theta$ implies that~$a_\theta=0$.
Thus, if~$b$ vanishes identically, then~$a$ is constant.  
As will be seen in \S\ref{sssec: lochomexamp}, 
this implies that~$M\subset X$ is locally homogeneous.

\textsl{Invariants in normal coordinates.}
Another interpretation of the invariants~$a$ and~$b$ (and the
canonical coframing) can be seen by considering a defining function
for the hypersurface in so-called \emph{normal coordinates}.

According to~\cite[\S4.2]{MR2000b:32066}, 
a holomorphic coordinate system~$(w,z)$ centered on~$p\in M\subset X$
is said to be \emph{normal} if there is an open $p$-neighborhood~$U$
in~$X$ such that $M\cap U$ is defined by an equation of the form
\be\label{eq: normalcoords}
\text{Im}\,w = z\bar z\,\psi\bigl(z,\,\bar z,\,\text{Re}\,w\,\bigr).
\ee
for some smooth function~$\psi$ of its arguments.

Acording to~\cite[Theorem~4.2.6]{MR2000b:32066}, 
when~$M$ is real-analytic, normal coordinates centered on~$p$ exist
and the function~$\psi$ can be regarded as an analytic function
of its arguments.  If~$M$ is merely smooth, then, for any
integer~$k$, there exist~$p$-centered coordinates such
that~\eqref{eq: normalcoords} holds up to order~$k$.

There are two limitations of normal coordinates for our purposes:
First, the proofs given in~\cite{MR2000b:32066} do not show that
one can choose the coordinates to be unimodular and, second, normal
coordinates are not unique.  

It is not difficult to show, however, that, for any point~$p\in M$,
there exist $p$-centered unimodular holomorphic coordinates~$(w,z)$
with~$w = u + \iC\,v$ such that~$M$ is defined in these coordinates
by an equation of the form~\eqref{eq: definingF} where
\be
F = {\ts\frac12}\,z\bar z\,\bigl(1 + b(p)\,z^2 
                 + {\ts\frac32} a(p)\,z\bar z 
                 +  \overline{b(p)}\,{\bar z}^2 \bigr)
    + {\mathrm{O}}(5)
\ee
and where the unwritten terms vanishing to order~$5$ and higher 
may contain~$u$ as well as~$z$ and~$\bar z$.%
\footnote{In the real-analytic case, it is possible to 
canonically define unique $p$-centered unimodular coordinates
by requiring that they satisfy certain natural geometric
conditions.  However, these coordinates do not seem to be
particularly useful and their description is less than
satisfactory, so their definition has not been included here. }

\textsl{CR-invariants.}
Of course, the usual (i.e., nonunimodular) CR-invariants of~$M$ 
can be expressed in terms of unimodular CR-invariants.%
\footnote{For the reader's convenience, the definition of
the Cartan connection of a real hypersurface in a complex
surface and its invariants are reviewed 
in Section~\ref{sec: CRinvs}.}

\begin{proposition}[Criterion for CR-flatness]
\label{prop: CRflatcrit}
The {\upshape{CR}}-structure on~$M$ is {\upshape{CR}}-flat if and only if
\be
a_{\bar\eta\bar\eta} + 2\iC\,b_\theta + 6 a b = 0.
\ee
\end{proposition}

\begin{proof}
By the analysis conducted in Section~\ref{sec: CRinvs}, one knows
that there exist $1$-forms~$\alpha$, $\beta$, $\sigma = \overline\sigma$
and a function~$s$ on~$M$ such that the following equations hold
\be\label{eq: CRflattest}
\begin{aligned}
\d\theta &=-(\alpha+\bar\alpha)\w\theta + \iC\,\eta\w\bar\eta\,,\\
\d\eta &= -\beta\w\theta-\alpha\w\eta\,,\\
\d\alpha &= -\sigma\w\theta-\iC\,\beta\w\bar\eta-2\iC\,\bar\beta\w\eta\,,\\
\d\beta &= -\sigma\w\eta+\bar\alpha\w\beta - s\,\overline{\eta}\w\theta\,.\\
\end{aligned}
\ee
Moreover, the underlying CR-structure is CR-flat if and only if the
function~$s$ vanishes identically.

Thus, to verify the proposition, it suffices to determine 
a triple of forms $\alpha$, $\beta$, and~$\sigma = \overline{\sigma}$ 
and a function~$s$ satisfying~\eqref{eq: CRflattest} 
in terms of the primary invariants and their derivatives.

Now, \eqref{eq: CRflattest} does not uniquely determine the 
unknown forms $\alpha$, $\beta$ and~$\sigma = \overline{\sigma}$, 
but, following the calculations in Section~\ref{sec: CRinvs},
any solution is seen to be of the form
\be
\begin{aligned}
\alpha &= (u - {\ts\frac32}\iC\,a)\,\theta\,,\\
\beta &=  (u + {\ts\frac12}\iC a)\,\eta + 2\iC b\,\overline{\eta}
            -a_{\bar\eta}\,\theta\,,\\ 
\sigma &= -\d u 
-{\ts\frac12}\iC\bigl(a_\eta\,\eta - a_{\bar\eta}\,\overline{\eta}\bigr)
            +(u^2+{\ts\frac14}\,a^2- 4b\overline{b} 
              -{\ts\frac12} a_{\eta\bar\eta}
              -{\ts\frac12} a_{\bar\eta\eta})\,\theta\,,\\
s &= a_{\bar\eta\bar\eta} + 2\iC\,b_\theta + 6 a b\,.
\end{aligned}
\ee
for some real-valued function~$u$ on~$M$.  (In particular, taking~$u=0$
gives the only solution that abides~$\alpha = -\overline{\alpha}$.)
Consequently, the underlying CR-structure is CR-flat if and only if
$a_{\bar\eta\bar\eta} + 2\iC\,b_\theta + 6 a b=0$, as claimed.
\end{proof}

\subsubsection{Locally homogeneous examples}\label{sssec: lochomexamp}
A real hypersurface~$M\subset X$ is said to be locally homogeneous, 
if, for any two points~$p,q\in M$, there exists an open 
$p$-neighborhood~$U\subset X$ and a unimodular 
biholomorphism~$\phi:U\to \phi(U)\subset X$ such that~$\phi(p)=q$
and~$\phi(M\cap U) = M\cap \phi(U)$.

If~$M$ is locally homogeneous, then the canonical 
coframing~$(\theta,\eta)$ must also be locally homogeneous.  
Consequently, the primary invariants~$a$ and~$b$ 
must be constant functions on~$M$.

Conversely, for any constants~$a = \bar a$ and~$b$, there exists
a connected, simply connected $3$-dimensional Lie group~$G_{a,b}$
on which there exist left-invariant $1$-forms~$\theta=\bar\theta$
and~$\eta$ satisfying~$\iC\theta\w\eta\w\overline{\eta}\not=0$ and
\be
\d\theta = \iC\eta\w\overline{\eta}
\qquad\text{and}\qquad
\d\eta = 2\iC\,\theta\w\bigl( a\,\eta + b\,\overline{\eta}\,\bigr).
\ee
The group~$G_{a,b}$ and the coframing~$(\theta,\eta)$ are unique
up to isomorphism. 

Thus, there is a $3$-parameter family of homogeneous abstract
unimodular CR-structures.   Note that, by Proposition~\ref{prop: CRflatcrit},
the underlying CR-structure is CR-flat if and only if either~$a=0$
or~$b=0$.  

When~$a>|b|$, 
the group~$G_{a,b}$ is isomorphic to $\SU(2)$. 
When~$a = \pm|b|$ but~$(a,b)\not=(0,0)$, 
the group~$G_{a,b}$ is isomorphic
to a semi-direct product of~$\bbR$ with~$\bbR^2$.
The group~$G_{0,0}$ is isomorphic to the Heisenberg group.
In all other cases, the group~$G_{a,b}$ is isomorphic 
to the simply-connected cover of $\SL(2,\bbR)$.

\begin{example}[The case~$(a,b)=(0,0)$]
Let~$(w,z)$ be standard coordinates on~$\bbC^2$, 
let~$\Upsilon=\d w\w\d z$, and let~$w = u + \iC v$.
Then, for the hypersurface~$M\subset\bbC^2$
defined by~$v = \mathrm{Im}\,w = \frac12\,z\bar z$, 
one computes
\be
\theta = \d u +{\ts\frac{\iC}2}\,(z\,\d\bar z- \bar z\,\d z)
\qquad\text{and}\qquad
\eta = \d z,
\ee
so that the structure equations are
\be
\d\theta = \iC\,\eta\w\overline{\eta}
\qquad\text{and}\qquad
\d\eta = 0.
\ee
Thus, the primary invariants~$a$ and~$b$ both vanish.

This example is seen to be homogeneous by inspection.  
The $3$-parameter group of unimodular biholomorphic transformations
defined by
\be
f(w,z) = \bigl(w + \iC\,\overline{z_0}\, z
            + u_0 +{\ts\frac\iC2}\,z_0\overline{z_0}
              \,,\,\,z + z_0\bigr)
\ee
for arbitrary~$z_0\in\bbC$ and~$u_0\in\bbR$ preserves~$M$ 
and acts transitively on it.  
\end{example}

\begin{example}[Homogeneous examples with~$a>|b|$]
\label{ex: SO3exs}
To construct these homogeneous examples, whose symmetry
groups are isomorphic to either~$\SU(2)$ or~$\SO(3)$, 
explicitly as hypersurfaces in unimodular complex surfaces, 
one can exploit the geometry of~$\SO(3)$-actions.

Let~$\SO(3)$ act on~$\bbC^3$ as the complexification 
of its standard action on~$\bbR^3$.  This action leaves
invariant the holomorphic quadratic function~$Q(z)=z\cdot z$ 
and the positive definite Hermitian form~$H(z) = z\cdot\bar z
= |z|^2$.  These satisfy the obvious inequalities 
\be
|z|^2 \ge |z\cdot z|\ge 0.
\ee
The SO(3)-orbit of a nonzero~$z_0$ has dimension~$3$ 
if and only if $|z_0|^2 > |z_0\cdot z_0|$ since
$|z_0|^2 = |z_0\cdot z_0|>0$ 
if and only if~$z_0 = \lambda x_0$ for some~$\lambda\in\bbC^*$ 
and some nonzero~$x_0\in\bbR^3$.  In fact,~$\SO(3)$
acts freely on the open set where~$|z|^2 > |z\cdot z|$.

Thus, for constants~$h\in\bbR^+$ and~$q\in\bbC$ 
with~$h>|q|$, the submanifold~$M_{h,q}$ defined by
\be\label{eq: SO3invMhqdef}
M_{h,q} = \bigl\{ z\in\bbC^3 \mid z\cdot z = q,\ |z|^2 = h\ \bigr\}
\ee
is an $\SO(3)$-invariant, $3$-dimensional hypersurface 
in the $\SO(3)$-invariant quadric surface~$X_q\subset\bbC^3$
defined by~$z\cdot z = q$.  It is not difficult to see that~$M_{h,q}$
is CR-nondegenerate and is an~$\SO(3)$-orbit.

On the other hand, the locus~$M_{|q|,q}\subset X_q$ defined by
setting~$h=|q|$ in~\eqref{eq: SO3invMhqdef} is either a point
(in the case~$q=0$) or a smoothly embedded, totally real 
$2$-sphere (in the case~$q\not=0$).  These singular orbits will
become significant when the behavior of the homogeneous
hypersurfaces under the unimodular normal flow (yet to be 
defined) is studied.

Now, consider the $\SO(3)$-invariant $2$-form
\be
\Upsilon 
= \frac{z^1\,\d z^2\w\d z^3+z^2\,\d z^3\w\d z^1 +z^3\,\d z^1\w\d z^2}
  {(z^1)^2+(z^2)^2+(z^3)^2}\,,
\ee 
which is well-defined away from the singular quadric~$z\cdot z = 0$.
It is not difficult to see that~$\Upsilon$ pulls back to each
quadric surface~$X_q$ with~$q\not=0$ to be a nowhere-vanishing
$2$-form~$\Upsilon_q$, one that is invariant under the action 
of~$\SO(3)$.

Since $M_{h,q}$ is homogeneous under the unimodular action 
of~$\SO(3)$ on~$(X_q,\Upsilon_q)$,
its primary invariants are constant.  
Direct computation yields
\be\label{eq: abasrxi}
a = \frac{ h}{2\bigl(h^2 - |q|^2\bigr)^{2/3}}
\qquad\text{and}\qquad
b = \frac{-q}{2\bigl(h^2 - |q|^2\bigr)^{2/3}}.
\ee
These formulae can be solved for~$h$ and~$q$ in the form
\be\label{eq: hqasab}
h = \frac{ a}{8\bigl(a^2 - |b|^2\bigr)^2}
\qquad\text{and}\qquad
q = \frac{-b}{8\bigl(a^2 - |b|^2\bigr)^2}.
\ee
Consequently, any~$(a,b)\in\bbR\times\bbC$ with~$a>|b|>0$ 
is represented by a unique~$M_{h,q}$ with~$h>|q|>0$.

The reader may wonder what happens in the case~$q=0$, which is,
apparently, not covered by the above formulae.  It is an 
interesting fact that, even though~$\Upsilon$ is not well-defined 
on the locus~$z\cdot z = 0$,
there is a well-defined~$\Upsilon_0$ on~$X_0$ that is, in a
natural sense, the limit of~$\Upsilon_q$ on~$X_q$ as~$q\to 0$
and that is a nonvanishing $2$-form on the smooth part of~$X_0$,
i.e., away from~$z=0$.  One sees this by noticing that, on the open set
where~$z^1\not=0$, the $2$-form
\be
\Psi_1 = \frac{\d z^2\w\d z^3}{z^1}
\ee
pulls back to each level set~$X_q$ with~$q\not=0$ to agree with
$\Upsilon_q$.  Moreover,~$\Psi_1$ pulls back to be well-defined 
and nonvanishing in the open set in~$X_0$ on which~$z^1\not=0$ 
and thus can naturally be regarded as the limit~$\Upsilon_0$ 
of~$\Upsilon_q$ as~$q\to 0$ on this open set.  By cyclically 
permuting the indices, one sees that~$\Upsilon_0$ can be defined 
in this manner on all of~$X_0$ away from the singular point~$z = 0$.  
Taking this~$\Upsilon_0$ as the volume form on~$X_0$, the 
formulae~\eqref{eq: abasrxi} are seen to be valid even for~$q=0$.

One more comment about this construction:  For~$q\not=0$,
the surface~$X_q$ is smooth, connected, and simply-connected.
When~$h>|q|$, the hypersurface~$M_{h,q}\subset X_q$ 
is compact and bounds the compact, simply-connected domain 
defined by~$|z|^2\le h$, even though~$M_{h,q}$, being 
diffeomorphic to~$\SO(3)$, is not, itself, simply-connected. 
As will be seen below, it is very unlikely that
the simply-connected cover of~$M_{h,q}$ can be embedded
as a hypersurface in a unimodular surface in such a way that
it bounds a compact domain.

On the other hand, the surface~$X_0$ is singular at~$z=0$
and this singularity is what is known as a `double point'
in algebraic geometry.  In fact, one can realize~$(X_0,\Upsilon_0)$
as the quotient of~$(\bbC^2,\d w\w\d z)$ by the unimodular
involution~$(w,z)\to(-w,-z)$.  Under this identification,
the~$SO(3)$-action on~$X_0$ lifts to the standard~$\SU(2)$-action 
on~$\bbC^2$ (which is, of course, unimodular).  The level sets~$M_{h,0}$
are the quotients of the hyperspheres centered on~$(0,0)\in\bbC^2$.
Thus, the simply-connected models for the case
of constant invariants~$a>0$ and~$b=0$ can be embedded the standard 
hyperspheres in~$\bbC^2$.  For more detail on this case, see
Example~\ref{ex: concenspheres}.
\end{example}

\begin{example}[Homogeneous examples with~$a<0$ and~$b=0$]
\label{ex: homa<0b=0}
Let~$\bbD$ be the unit disk in~$\bbC$ endowed with
the Poincar\'e metric of constant curvature~$-1$,
whose oriented isometry group is~$\PSL(2,\bbR)$.
Let~$X=T^*\bbD\to\bbD$
be the canonical bundle endowed with the standard (nowhere-vanishing)
holomorphic $2$-form~$\Upsilon$ that it inherits as a (holomorphic)
cotangent bundle.  Let~$|\cdot|:X\to\bbR$ denote the hermitian norm
on~$X$ that it inherits from the Poincar\'e metric.

For~$r>0$, let~$M_r = \{\zeta\in X\ \vrule\ |\zeta|^2 = \frac12\,r^2\}$.
Then~$M_r$ is a hypersurface in~$X$ that is easily seen to be
nondegenerate and one computes without difficulty 
that the invariants satisfy~$a = - \frac12 r^{-2/3}$ and~$b=0$.

Of course,~$M_r$ is not simply-connected but is diffeomorphic 
to~$S^1\times\bbD$.  One can get the simply-connected example
by considering the inverse image of~$M_r$ under the universal
covering~$\tilde X^*\to X^*$ where~$X^*\subset X$ is the
complement of the zero section of~$X\to\bbD$.

In the other direction, for any group~$\Gamma\subset\PSL(2,\bbR)$ 
acting freely on~$\bbD$ with compact quotient~$\Gamma\backslash\bbD$,
the action extends unimodularly and freely to~$(X,\Upsilon)$. 
The quotient~$\Gamma\backslash M_r$ is then a compact example
of a nondegenerate hypersurface in~$(\Gamma\backslash X,\Upsilon)$
that has constant invariants and that bounds a compact domain. 

For more discussion about this example, 
see Example~\ref{ex: homa<0b=0evol}.
\end{example}

\begin{example}[Homogeneous examples with~$a<-|b|$ or~$|a|<|b|$]
\label{ex: PSL2exs}
Examples of these two types have symmetry group isomorphic to 
some covering of~$\mathrm{PSL}(2,\bbR)=\SL(2,\bbR)/\{\pm\mathrm{I}_2\}$.
To construct them up to a covering, 
one can exploit the geometry of the irreducible
$3$-dimensional representation of~$\PSL(2,\bbR)$.

Let~$\SL(2,\bbR)$ act on the space~$S_2(\bbR)$ of $2$-by-$2$ symmetric
matrices with real entries in the usual way:  $A\cdot s = A\,s\,{}^t\!A$
for~$A\in\SL(2,\bbR)$ and~$s\in S_2(\bbR)$.  This action is almost
faithful, with~$\{\pm\mathrm{I}_2\}\simeq\bbZ_2$ acting trivially,
so that this is actually a representation of~$\PSL(2,\bbR)$.
This action preserves the quadratic form~$s\mapsto \det(s)$, which has
maximal negative definite subspaces of dimension~$2$.  Let~$r\cdot s$
denote the inner product on~$S_2(\bbR)$ associated to this quadratic
form, i.e.,~$s\cdot s = \det (s)$.  Because~$\PSL(2,\bbR)$ is
simply connected, it also preserves an orientation of~$S_2(\bbR)$
as well as the `positive time-like cone' $S^+_2(\bbR)$ that consists 
of the positive definite matrices.

Let~$\PSL(2,\bbR)$ act on~$S_2(\bbC)=\bbC\otimes S_2(\bbR)$ by 
extension of scalars and extend the inner product on~$S_2(\bbR)$ 
complex linearly to~$S_2(\bbC)$.  Then~$\PSL(2,\bbR)$ preserves
both the complex quadratic form~$Q(z) = z\cdot z$ and the 
Hermitian form~$H(z) = z\cdot \bar z$.  For~$q\in\bbC$, define
\be
X_q = Q^{-1}(q) = \{\,z\in S_2(\bbC)\mid z\cdot z = q\,\},
\ee
so that~$X_q$ is an $\SL(2,\bbR)$-invariant quadric surface,
nonsingular if~$q\not=0$ while~$X_0$ is a (singular) cone.

Just as in the $\SO(3)$-invariant case that has already been
treated, one can define a canonical, $\PSL(2,\bbR)$-invariant
holomorphic $2$-form~$\Upsilon_q$ on~$X_q$ so that, when~$q\not=0$,
the inclusion mapping~$z:X_q\to S_2(\bbC)\simeq\bbC^3$ satisfies
\be
{\ts\frac12}\,z\w \d z\w \d z = q \Upsilon_q\,\mathbf{v}
\ee
where~$\mathbf{v}\in \Lambda^3\bigl(S_2(\bbR)\bigr)$ is
the volume $3$-vector associated to an orthnormal, oriented
basis of~$S_2(\bbR)\simeq\bbR^3$ and so that, when~$q=0$, the
$2$-form~$\Upsilon_0$, well-defined on the smooth part of~$X_0$,
is the limit of~$\Upsilon_q$ as~$q$ approaches~$0$ in a suitable
sense.  (The details, which are entirely analogous to those in
the $\SO(3)$ case, will be left to the reader).  
In this way,~$(X_q,\Upsilon_q)$ becomes a $\PSL(2,\bbR)$-invariant%
\footnote{In fact, $\PSL(2,\bbC)$-invariant.} 
unimodular complex surface.

Now, the $\PSL(2,\bbR)$ action must preserve the
level sets of~$H$: For~$h\in\bbR$, set
\be
M_{h,q} = Q^{-1}(q)\cap H^{-1}(h) 
        = \{\,z\in X_q\mid z\cdot \bar z = h\,\}.
\ee

The triangle inequality for the indefinite quadratic form~$\det$
on~$S_2(\bbR)$ implies that $H(z) \le |Q(z)|$.  In other words,
$M_{h,q}$ is empty unless~$h\le |q|$.  It is sometimes useful 
to note the identity
\be\label{eq: Mhqscaling}
\lambda\cdot M_{h,q} = M_{|\lambda|^2 h,\,\lambda^2 q}\,,
\ee
valid for~$\lambda\in\bbC^*$,
and the fact that scalar multiplication by~$\lambda$ 
induces a mapping from~$X_q$ to~$X_{\lambda^2 q}$
that pulls back~$\Upsilon_{\lambda^2 q}$ to~$\lambda\,\Upsilon_q$.
 
Computation shows that each~$M_{h,q}$ with~$h<|q|$ 
and~$h\not=-|q|$ is a (nonempty) smooth, nondegenerate 
real hypersurface in~$(X_q,\Upsilon_q)$ 
that has constant invariants
\be\label{eq: abashq}
a = \frac{h}{2(h^2-|q|^2)^{2/3}}
\qquad\text{and}\qquad
b = \frac{-q}{2(h^2-|q|^2)^{2/3}}\,.
\ee
Referring to the inversion formulae~\eqref{eq: hqasab},
one sees that, when~$a<-|b|$, equations~\eqref{eq: abashq} 
can be solved for~$h$ and~$q$ uniquely with~$h < -|q|$.  On the
other hand, when~$|a|<|b|$, the equations~\eqref{eq: abashq} 
can be solved for~$h$ and~$q$ uniquely with~$|h| < |q|$.
Thus, these provide homogeneous (though not 
simply connected) models of the desired kinds.

Now, when~$h^2\not=|q|^2$, the hypersurface~$M_{h,q}$ 
is not connected.  

When~$h<-|q|$, any element
of~$M_{h,q}$ is of the form~$x-\iC y$ where~$x$ and~$y$
in~$S_2(\bbR)$ are linearly independent and span a negative
definite $2$-dimensional subspace of~$S_2(\bbR)$.  Say
that~$x-\iC y$ lies in~$M^+_{h,q}$ if the basis~$( \mathrm{I}_2,x,y)$
is a positively oriented basis of~$S_2(\bbR)$ and 
that $x-\iC y$ lies in~$M^-_{h,q}$ if the basis~$(-\mathrm{I}_2,x,y)$
is a positively oriented basis of~$S_2(\bbR)$.  Then~$M_{h,q}$
is the disjoint union of the connected sets~$M^\pm_{h,q}$.
In fact, $\PSL(2,\bbR)$ acts simply transitively 
on each of~$M^+_{h,q}$ and~$M^-_{h,q}$, as is not difficult
to see. 

On the other hand, when~$|h|<|q|$, the situation is a bit 
more subtle:  The hypersurface~$M_{h,q}$ has two components,
but there is no natural way to label one as positive and
the other as negative.  This is most easily seen as follows:
Consider the case~$q=1$.  Each element of~$M_{h,1}$ is of
the form~$x-\iC y$ where~$x,y\in S_2(\bbR)$ satisfy 
\be
x\cdot x = {\ts\frac12}(h+1) > 0, \quad
y\cdot y = {\ts\frac12}(h-1) < 0, \quad\text{and}\quad
x\cdot y =  0.
\ee 
Now, one can define
$M^+_{h,1}$ to be those elements~$x-\iC y\in M_{h,1}$ 
such that~$x$ lies in~$S^+_2(\bbR)$ and 
$M^-_{h,1}$ to be those elements~$x-\iC y\in M_{h,1}$ 
such that~$-x$ lies in~$S^+_2(\bbR)$.  Then~$M_{h,1}$
is the disjoint union of the two connected components
$M^+_{h,1}$ and~$M^-_{h,q}$.  Moreover, $\PSL(2,\bbR)$ 
acts simply transitively on each of~$M^+_{h,1}$ and~$M^-_{h,1}$.
Now, using the scaling property~\eqref{eq: Mhqscaling},
one might be tempted to try to `define'~$M^\pm_{h,q}$
as the image under scalar multiplication by~$\sqrt{q}\in\bbC^*$
of~$M^\pm_{h/|q|,1}$.  However, the choice of the square root
affects which component is mapped to which component, and,
of course, there is no continuous way to choose this square root. 

For use in the study of the effect of unimodular normal
flow on these examples, it is also worthwhile to look at 
the sets~$M_{h,q}$ when~$h = \pm|q|$.   

The case~$M_{0,0}$ is special.  It is not difficult to
see that
\be
M_{0,0} 
= \left\{\ \lambda\,x \mid \lambda\in\bbC,\ x\in S_2(\bbR)
\ \text{and}\ \det(x)=0\ \right\}.
\ee
Thus,~$M_{0,0}$ is a union of a real $1$-parameter family of 
complex lines.  In particular, $M_{0,0}$ is a totally degenerate
real hypersurface in~$X_0$.  It is not homogeneous under~$\PSL(2,\bbR)$,
as the orbits are $2$-dimensional, except for the origin itself,
which is $0$-dimensional.

The subset~$M_{1,1}$ is  described as
\be
M_{1,1} = \left\{\ x\in S_2(\bbR)\mid\det(x)=1\ \right\},
\ee
and hence is a $2$-dimensional surface, the (real) unit hyperboloid
of $2$-sheets, sometimes denoted by~$H_1$.  
Its two sheets are~$M^+_{1,1}=H^+_1$, 
consisting of the positive definite~$s\in S_2(\bbR)$ 
with determinant~$1$, and~$M^-_{1,1} = - M^+_{1,1} = - H^+_1 = H^-_1$.

The subset~$M_{-1,1}$ is somewhat more complicated:  It is
$3$-dimensional and singular.  It can be described in equations
in the form
\be
M_{-1,1}
 = \left\{\ x-\iC\,y \mid
   x,y\in S_2(\bbR),\ x{\cdot}x=x{\cdot}y=0,\,y{\cdot}y = -1\ \right\}.
\ee
Thus,~$M_{-1,1}$ is the union of two smooth, closed, connected 
$3$-dimensional hypersurfaces in~$X_1$, say~$M^{\pm}_{-1,1}$, each of which 
is seen to be an oriented line bundle over the hyperboloid of one sheet
\be
H_{-1} = \left\{\ y\in S_2(\bbR) \mid \det(y)=-1\ \right\}
\ee
via the mapping~$x-\iC\,y \mapsto y$.  (Each fiber of the map
$M_{-1,1}\to H_{-1}$ is a union of two distinct null lines in~$S_2(\bbR)$
and these lines are oriented by taking the positive ray to lie
in the closure of the positive cone~$S^+_2(\bbR)$.)

The (CR-degenerate) hypersurfaces~$M^+_{-1,1}$ and~$M^-_{-1,1}$
intersect transversely along the singular locus of~$M_{-1,1}$, 
\be
\mathrm{sing}\bigl(M_{-1,1}\bigr)
 = \left\{\ \iC\,y\in S_2(\bbR)\ 
         \mid\ \det(y)=-1\ \right\} = \iC\, H_{-1}.
\ee
The smooth part~$M^*_{-1,1}$ of~$M_{-1,1}$ is the disjoint
union of $4$ components, which can be labeled as~$M^{\pm,\pm}_{-1,1}$.
Each of these components is acted upon simply transitively
by~$\PSL(2,\bbR)$.
\end{example}

\subsubsection{Cohomogeneity one examples}\label{sssec: cohom1exs}
The homogeneous examples display many interesting
phenomena, but they are necessarily of somewhat limited use
in understanding the general case.  A somewhat more interesting
case is that of examples that have cohomogeneity one, i.e., 
the group of symmetries acts with principal orbits of dimension~$2$.
It turns out that, in this case, the structure equations can
be analyzed fully and these examples characterized.

If the pseudogroup of local symmetries of~$(M,\Phi)$ has
$2$-dimensional orbits on~$M$, then, using~\eqref{eq: deta}, 
the invariants of the canonical coframing satisfy
\be\label{eq: dadbmatrix}
\bpm \d a\\ \d b\\ \d \bar b\epm
= \bpm 
a_\theta& a_\eta& \overline{a_\eta}\\
b_\theta& \overline{a_\eta}& b_{\bar\eta}\\
\overline{b_\theta}& \overline{b_{\bar\eta}}& a_\eta
\epm
\bpm \theta\\ \eta\\ \bar\eta\epm,
\ee
where the $3$-by-$3$ matrix of covariant derivatives
in~\eqref{eq: dadbmatrix} must have rank at most~$1$
since~$(a,b):M\to\bbR\times\bbC$ must be constant on
the $2$-dimensional orbits of the pseudogroup of symmetries.

Now, the rank assumption implies, in particular,
that~$|a_\eta|^2 = |b_{\bar\eta}|^2$.  
Thus, if the locus~$a_\eta=0$ has any interior, it will
follow that~$\d a$ and $\d b$ are multiples of~$\theta$
on that interior and hence, since~$\theta$ is a contact form,
that~$\d a$ and~$\d b$ must vanish on the interior of 
the locus~$a_\eta=0$.  In particular,~$a$
and~$b$ must be constant on each component of the
interior of the locus~$a_\eta=0$.  Consequently, the
unimodular CR-structure will be locally homogeneous on
each such component.  For simplicity, this case will be
set aside and it will be assumed for the rest of this
subsubsection that the locus~$a_\eta=0$ has no interior.

Let~$M^*\subset M$ be the open set on which~$a_\eta\not=0$
and write
\be
\d a = r\left(\,2s\,\theta + \E^{\iC\phi}\,\eta 
             + \E^{-\iC\phi}\,\overline{\eta}\,\right)
\ee
for some real functions~$r = |a_\eta|>0$,~$s$, and
some function~$\phi$ defined up to an integral multiple
of~$2\pi$ on~$M^*$.   Since~$r$, $s$, and~$\phi$ are 
also invariant under the action of the symmetry pseudogroup,
their differentials must be multiples of~$\d a$.  
In particular, the $1$-form
\be\label{eq: cohom1alphadef}
\alpha = 2s\,\theta + \E^{\iC\phi}\,\eta 
             + \E^{-\iC\phi}\,\overline{\eta}
\ee
must be closed.  Moreover, since~$\d b$ must be a multiple
of~$\d a$, and, in particular,~$a_{\bar\eta}\,\d a - a_\eta\,\d b = 0$,
it follows that
\be
\d b = \E^{-2\iC\phi}\,\d a = r\E^{-2\iC\phi}\,\alpha.
\ee

Now, computing the exterior derivatives of the structure equations
and the identity~$\d(\d\alpha)=0$ yields that, first
\be
\d\phi = s\,\alpha,
\ee
and then, that
\be
b = \E^{-2\iC\phi}\left( a + s^2 + \iC\,v\right)
\ee
for some real-valued function~$v$.  Moreover, one finds that
\be
\d s = -v\,\alpha
\qquad\text{and}\qquad
\d v = 2s(a + s^2)\,\alpha.
\ee

In conclusion, one has structure equations of the form
\be\label{eq: CartanStrEqsCohom1}
\begin{aligned}
\d\theta &= \iC\,\eta\w\overline{\eta},\\
\d\eta &= 2\iC\,\theta\w\bigl(a\,\eta 
 + \E^{-2\iC\phi}\left( a + s^2 + \iC\,v\right)\,\overline{\eta}\bigr),\\
\bpm\d a\\ \d v\\ \d s\\ \d\phi\epm
&=\bpm r\\ 2s(a{+}s^2)\\ -v\\ s \epm
\bigl(\,2s\,\theta + \E^{\iC\phi}\,\eta 
             + \E^{-\iC\phi}\,\overline{\eta}\,\bigr).
\end{aligned}
\ee
Note that the exterior derivatives of these
structure equations imply
\be
\d\bigl(\,2s\,\theta + \E^{\iC\phi}\,\eta 
             + \E^{-\iC\phi}\,\overline{\eta}\,\bigr)=
\d r \w \bigl(\,2s\,\theta + \E^{\iC\phi}\,\eta 
             + \E^{-\iC\phi}\,\overline{\eta}\,\bigr) = 0,
\ee
but no other identities on the derivatives of the forms~$\theta$
and~$\eta$ or the functions~$a$, $v$, $s$, or~$\phi$.
According to Cartan's generalization of the third fundamental
theorem of Lie, the solutions of~\eqref{eq: CartanStrEqsCohom1} 
depend on one arbitrary function of one variable up to diffeomorphism
and any coframing~$(\theta,\eta)$ satisfying these structure
equations has its pseudogroup acting transitively on the
leaves of the integrable $1$-form~$\alpha$ as defined 
in~\eqref{eq: cohom1alphadef}.

\subsection{Embedding results}\label{ssec: embeddingresults}
It is natural to ask how much of the local geometry of~$M\subset X$
is captured by the coframing~$(\theta,\eta)$.  In this section, 
some of these questions will be answered, at least in the 
real-analytic category.  First, it will be useful to have
an abstract notion of unimodular CR-structure:

\begin{definition}[Unimodular CR-structures]
\label{def: unimodCRstrs}
A closed, complex-valued $2$-form~$\Psi$ on a $3$-manifold~$M$
will be said to be a \emph{unimodular CR-structure} on~$M$ 
if the real and imaginary parts of~$\Psi$ are linearly independent
at each point of~$M$.  Given such a~$\Psi$, 
the unique $2$-plane field~$D\subset TM$ that consists of the
$2$-planes~$D_x\subset T_xM$ to which~$\Psi$ pulls back to zero
will be called the \emph{complex tangent space} of~$\Psi$.  
If~$D$ is a contact plane field, then $\Psi$ will be said to
be \emph{nondegenerate}.  
\end{definition}

\begin{remark}[Canonical coframing for unimodular CR-structures]
It is clear that, if~$\Psi$ is a nondegenerate unimoduclar CR-structure,
then the method employed in the proof of Proposition~\ref{prop: cancoframing}
will prove that there exists a unique coframing~$(\theta,\eta)$
with~$\theta=\bar \theta$ such that~$\Psi = \theta\w\eta$ and
$\d\theta = \iC\,\eta\w\overline{\eta}$.  This coframing will, of
course, be referred to as the canonical coframing of~$(M,\Psi)$.
\end{remark}

The first result provides local uniqueness for realizing
a real-analytic unimodular CR-structure 
as a hypersurface in a unimodular complex surface:

\begin{proposition}[Ambient uniqueness]\label{prop: ambunique}
Let~$M^3$ be a real-analytic $3$-manifold endowed with a
real-analytic unimodular CR-structure~$\Psi$.

Suppose that~$(X_1,\Upsilon_1)$ and~$(X_2,\Upsilon_2)$ are
two unimodular complex surfaces such that there exist 
real-analytic embeddings~$\phi_i:M\to X_i$ satisfying
\be
{\phi_1}^*(\Upsilon_1) = {\phi_2}^*(\Upsilon_2) = \Psi.
\ee
Then there exist open sets~$U_i\subset X_i$, 
with~$U_i$ containing $\phi_i(M)$ and a unimodular
biholomorphism~$\Phi:(U_1,\Upsilon_1)\to(U_2,\Upsilon_2)$ 
such that~$\phi_2 = \Phi\circ\phi_1$.

Moreover, $(\Phi,U_1,U_2)$ is locally unique in the 
sense that, if~$(\tilde\Phi,\tilde U_1,\tilde U_2)$
is another triple with these properties, then there
exists an open subset~$V_1\subset U_1\cap \tilde U_1$
containing~$\phi_1(M)$ such that~$\tilde\Phi = \Phi$ on~$V_1$.
\end{proposition}

\begin{proof}
The proof is an application of the Cartan-K\"ahler Theorem:
Consider the product manifold~$X_1\times X_2$ 
with projections~$\pi_i: X_1\times X_2\to X_i$ 
and the ideal~$\cI$ on~$X_1\times X_2$ generated by the 
real and imaginary components 
of the closed~$2$-form~$\pi_2^*\Upsilon_2 - \pi_1^*\Upsilon_1$.
The ideal~$\cI$ 
with independence condition~$\Upsilon_1\w\overline{\Upsilon_1}\not=0$ 
is evidently involutive, with Cartan
characters
\be
(s_0,s_1,s_2,s_3,s_4) = (0,2,2,0,0).
\ee

The graph~$\phi = (\phi_1,\phi_2):M\to X_1\times X_2$ is, 
by assumption, an integral manifold of~$\cI$ of dimension~$3$.
Inspection shows that it is regular and that, moreover, its
index of indeterminacy is zero.  By the Cartan-K\"ahler
theorem, it lies in a (locally unique) $4$-dimensional
integral manifold~$Y\subset X_1\times X_2$ of~$\cI$ to
which the $4$-form $\Upsilon_1\w\overline{\Upsilon_1}$
pulls back to be a volume form. 
It follows that, near~$\phi(M)$, 
the submanifold~$Y$ is the graph of a real-analytic invertible
mapping
$\Phi:U_1\to U_2$ (where~$U_i\subset X_i$ is an open neighborhood
of~$\phi_i(M)$) that satisfies~$\Phi^*\Upsilon_2 = \Upsilon_1$.  
Thus~$\Phi$ is a unimodular biholomorphism.  

The stated local uniqueness 
follows from the local uniqueness of~$Y$.
\end{proof}

\begin{proposition}[Ambient existence]\label{prop: ambexist}
Let~$M^3$ be a real-analytic $3$-manifold endowed with a
real-analytic unimodular CR-structure~$\Psi$.

Then there exists a unimodular complex surface~$(X,\Upsilon)$
and a real-analytic embedding~$\phi:M\to X$ 
such that~$\phi^*\Upsilon =\Psi$.
\end{proposition}

\begin{proof}
First, an application of the Cartan-K\"ahler Theorem will
show local existence:  On the product manifold~$M\times\bbC^2$,
consider the ideal~$\cI$ generated by the real and imaginary parts 
of the closed $2$-form~$\pi_2^*(\d z^1\w \d z^2)-\pi_1^*(\Psi)$.  
If~$\Omega$ is any volume form on~$M$, then the ideal~$\cI$ 
with independence condition~$\pi_1^*(\Omega)\not=0$ 
is involutive, with characters
\be
(s_0,s_1,s_2,s_3) = (0,2,2,0).
\ee

It follows that every point~$p\in M$ has an open neighborhood~$U$ 
on which there exists a real-analytic embedding~$z:U\to\bbC^2$
such that~$z^*(\d z^1\w\d z^2) = \Psi$.  

Second, a patching argument finishes the proof:
Consider the collection~$\mathcal{U}$ consisting of those open
sets~$U\subset M$ for which there exists a unimodular 
complex surface~$(X,\Upsilon)$ and a real-analytic 
embedding~$\phi:U\to X$ satisfying~$\phi^*\Upsilon = \Psi$.
By the first part of the argument, $\mathcal{U}$ is an open
covering of~$M$.  Let~$\mathcal{U}'$ be a countable, locally
finite refinement of~$\mathcal{U}$ and use the uniqueness
guaranteed by Proposition~\ref{prop: ambunique} to show that
the union of the elements of~$\mathcal{U}'$ belongs to~$\mathcal{U}$.
\end{proof}

\begin{remark}[Cohomogeneity 0 and 1 examples]
Suppose now that~$(M,\Psi)$ is a real-analytic, nondegenerate
unimodular CR-structure that is either of cohomogeneity $0$ 
(i.e., symmetry group~$G$ of dimension~$3$) 
or~$1$ (symmetry group~$G$ dimension~$2$).  It then follows from
Proposition~\ref{prop: ambunique} that 
the symmetry group action extends to any unimodular surface
thickening~$(X,\Upsilon)$ of~$(M,\Psi)$, at least to a neighborhood
of~$M\subset X$. Since~$M$ is a nondegenerate CR-structure,
the orbits of~$G$ cannot be complex curves and hence it follows
that the complexification~$G^\bbC$ of~$G$ extends to an infinitesimally
homogeneous unimodular action on~$(X,\Upsilon)$.  Thus, in the
usual way,~$(X,\Upsilon)$ can be regarded (up to a covering) 
as a coadjoint orbit of~$G^\bbC$.

In particular, the examples of cohomogeneity~$0$ and~$1$ can be
studied in a standard way by considering a $2$-dimensional
coadjoint orbit~$(X,\Upsilon)$ of the complexification of a 
$2$- or $3$-dimensional Lie group~$G$ and looking at the 
$G$-invariant hypersurfaces in these coadjoint orbits.

This gives a uniform way of interpreting the formulae found
in the previous subsection for the homogeneous examples.

Also, in the case that~$G$ is a $2$-dimensional Lie group
acting freely on~$(X,\Upsilon)$, the quotient~$G\backslash X$
will be a $2$-dimensional surface and the image of a $G$-invariant
hypersurface in~$X$ in this quotient will be a curve in this 
surface.  This gives a geometric interpretation to the result,
found by Cartan-K\"ahler analysis in~\S\ref{sssec: cohom1exs}, that the
cohomogeneity~$1$ examples depend on one function of one variable.

This point of view also simplifies the study of evolution of
such cohomogeneity~$1$ structures, as defined in the next section.
\end{remark}

\subsection[Flows]{Hypersurface flows in unimodular surfaces}
\label{ssec: flows}

In this section the behavior of the canonical coframing
and its invariants under deformation will be investigated.
A canonical flow for nondegenerate hypersurfaces will be
defined and some examples computed.

\subsubsection{Deformation formulae}
Consider a $1$-parameter family~$\iota:(-\epsilon,\epsilon)\times M\to X$ 
of immersions of~$M$ into~$X$ as nondegenerate hypersurfaces, 
parametrized so that
\be
\iota_*\left(\frac{\partial\hfil}{\partial t}(t,p)\right)
= f(t,p)\,J\bigl(T(t,p)\bigr),
\ee
for~$(t,p)\in (-\epsilon,\epsilon)\times M$, i.e., 
so that the variation is `normal'.  
The function~$f:(-\epsilon,\epsilon)\times M\to\bbR$ 
is the \emph{normal displacement function}.

By assumption,
\be
\iota^*(\Upsilon) = (\theta + \iC f\,\d t) \w \eta
\ee
where~$(\theta,\eta)$ is the ($t$-dependent) induced canonical coframing
on~$M$.  Expansion of the identities~$\d\bigl((\theta+\iC f\,\d t)\w\eta\bigr)
=\d(\d\theta)=0$ yields equations of the form
\be\label{eq: coframingevolf}
\begin{aligned}
\d\theta 
&= \iC\,\eta\w\overline{\eta} 
          +\bigl( {\ts\frac13}(4af+f_{\eta\bar\eta}+f_{\bar\eta\eta})\,\theta 
         + \iC(\,f_\eta\,\eta-f_{\bar\eta}\,\overline{\eta}\,)\bigr) \w \d t\\
\d\eta
&= 2\iC\,(\theta+ \iC\,f\,\d t)\w(a\,\eta + b\,\overline{\eta})
    +\bigl(\iC u\,\theta 
      - ({\ts\frac13}(4af+f_{\eta\bar\eta}+f_{\bar\eta\eta})
           +\iC f_\theta)\eta\bigr)\w\d t
\end{aligned}
\ee
where
\be
u = {\ts\frac13}(4af+f_{\eta\bar\eta}+f_{\bar\eta\eta})_{\overline{\eta}} 
        +\iC\ f_{\theta\overline{\eta}} + 4 f_\eta \overline{b}\,.
\ee

The evolution of the invariants~$a$ and~$b$ can be computed by
taking the exterior derivatives of the equations~\eqref{eq: coframingevolf}.

Note, in particular, that the normal deformation preserves 
the contact form~$\theta$ up to a multiple 
if and only if~$f_\eta = f_{\bar\eta}= 0$.
By~\eqref{eq: cofderexchformulae}, 
this implies that~$f_\theta=0$ as well, 
i.e., that~$f$ is a function of~$t$ alone.  In such a case, 
assuming that~$f(0)\not=0$, one can reparametrize in~$t$ 
so as to arrange that~$f(t)\equiv 1$ for~$t$ small.

\subsubsection{Unimodular normal flow}
The above formulae motivate the study of the case~$f\equiv 1$ for
a deformation.

\begin{definition}\label{def: unimodularnormalflow}
A $1$-parameter family~$\iota:(-\epsilon,\epsilon)\times M\to X$ 
of immersions of~$M$ into~$X$ as nondegenerate hypersurfaces will
be said to be a \emph{unimodular normal flow} if its normal
displacement function~$f$ satisfies~$f\equiv1$.
\end{definition}

For a unimodular normal flow, after normal reparametrization, one
has structure equations of the form
\be\label{eq: thetaetaunfformulae}
\begin{aligned}
\d\theta &= \iC\,\eta\w\overline{\eta} 
          +{\ts\frac43}a\,\theta\w \d t\\
\d\eta &= 2\iC\,(\theta+ \iC\,\d t)\w(a\,\eta + b\,\overline{\eta})
            +{\ts\frac43}\bigl(\iC a_{\bar\eta}\,\theta-a\,\eta\bigr)\w\d t,
\end{aligned}
\ee
while the primary invariants satisfy evolution equations of the form
\be\label{eq: abunftders}
\begin{aligned}
a_t &= {\ts\frac13}(a_{\eta\bar\eta}+a_{\bar\eta \eta}) 
          + {\ts\frac43}\,a^2 + 4\,b\overline{b}\\
b_t &= {\ts\frac23}\,a_{\bar\eta\bar\eta} + \iC\,b_\theta 
             + {\ts\frac{16}3}\,a b.
\end{aligned}
\ee

\begin{remark}[real-analytic uniqueness]
By~\eqref{eq: thetaetaunfformulae} and~\eqref{eq: abunftders},
for a unimodular normal flow, the first derivatives
of the coframing and primary invariants with respect to~$t$
are expressed in terms of the covariant derivatives of the invariants.
By induction, it follows that all higher time derivatives of
the coframing and primary invariants can be expressed in terms of
the coframing and iterated covariant derivatives of the invariants.
In particular, a real-analytic `initial' hypersurface can be embedded
into a real-analytic unimodular normal flow in at most one way.
Presumably, this is a special case of a uniqueness result for 
the flow in the smooth case.
\end{remark}

\begin{remark}[In unimodular coordinates]
Since the Reeb vector field~$T$ depends on three derivatives of
the embedding, the normal vector field~$\iC T$ also depends on 
three derivatives of the embedding.  In particular, the equation
for unimodular normal flow is, strictly speaking, a third order
flow.  

However, a closer examination reveals that, although~$\iC T$
is indeed a third order expression and is the lowest order 
normal vector field invariant under unimodular biholomorphism,
one can express~$\iC T$ (noninvariantly) as the sum of a
second order vector field and a vector field tangential to~$M$.
In other words, modulo reparametrization of~$M$ (in a possibly
non-normal fashion), one can write the flow as a second order
flow.

For example, suppose that the variation can be written locally
in unimodular coordinates~$(w,z)$ in the graphical form 
\be
v = F(z,\bar z, u, t)
\ee
where~$w = u + \iC\,v$.  Then one finds that this family represents 
a unimodular normal flow (after normal reparametrization) if and
only if~$F$ satisfies the nonlinear equation
\be\label{eq: evoleqngraphical}
F_t = \Bigl[2\left[
         (1{+}{F_u}^2)F_{z\bar z} + F_zF_{\bar z}F_{uu}
         - \iC (1{-}\iC\,F_u)F_{\bar z}F_{uz}
         + \iC (1{+}\iC\,F_u)F_{z}F_{u\bar z}
                    \right]\Bigr]^{1/3}.
\ee
Note that the expression inside the cube root is nonvanishing
as long as the family consists of nondegenerate hypersurfaces.

The linearization of~\eqref{eq: evoleqngraphical} at a nondegenerate
hypersurface~$v = F(z,\bar z, u)$ is an equation of the form
\be
G_t = S_F(G)
\ee
where~$S_F$ is a second-order subelliptic operator of the same 
symbol type as the $\bar\p_b$-Laplacian.  Thus, the linearization
is only weakly parabolic, so questions of short-time existence
or uniqueness do not appear to be easy to resolve.
\end{remark}

\begin{example}[Concentric spheres]\label{ex: concenspheres}
Consider the hypersurface foliation on~$\bbC^2$ minus the origin
defined by the level sets of the function~$r= |z|^2+|w|^2$, i.e.,
the spheres of positive radius centered on the origin.

Set
\be
\theta + \iC\,\d t = -\iC\,\frac{\bar z\,\d z + \bar w\,\d w}
                           {\bigl(|z|^2+|w|^2\bigr)^{1/3}}
\qquad\text{and}\qquad
\eta  =  \iC\,\frac{ z\,\d w - w\,\d z}
              {\bigl(|z|^2+|w|^2\bigr)^{2/3}}.
\ee
Then $(\theta + \iC\,\d t)\w\eta = \d z \w \d w$ and
\be
\d t = -\frac{\bar z\,\d z + \bar w\,\d w + z\,\d\bar z + w\,\d\bar w}
                           {2\bigl(|z|^2+|w|^2\bigr)^{1/3}}
     = -{\ts\frac12} r^{-1/3}\,\d r = \d\bigl(-{\ts\frac34}\,r^{2/3}\bigr).
\ee
Thus, the level sets of~$r>0$ are the integral surfaces of~$\d t$
and one sees that this is a unimodular normal flow when one sets
$t - t_0 = {\ts\frac34}{r_0}^{2/3}-{\ts\frac34}{r}^{2/3}$.

Computation yields
\be
\d\theta = \iC\,\eta\w\overline{\eta} + {\ts\frac43} r^{-2/3}\,\theta\w\tau,
\ee
so $(\theta,\eta)$ pulls back to each level set of~$r$ to be the 
canonical coframing on that level set.  The structure equation
for~$\eta$ turns out to be
\be
\d\eta = 2\iC\,\theta\w\bigl(r^{-2/3}\,\eta\bigr) 
         +{\ts\frac23}r^{-2/3}\,\eta\w\tau.
\ee

Note that the time~$t$ flow of the level set~$r=r_0>0$
is the level set~$r = R(t)$, where
\be
R(t) = \left({r_0}^{2/3} - {\ts\frac43}\,t\,\right)^{3/2},
\ee
and hence that this flow contracts to the center of the sphere
in finite time.
\end{example}

\begin{example}[Abstract Homogeneous Flows]\label{ex: abshomogflows}
Suppose now that~$M^3$ is endowed with an abstract
homogeneous unimodular CR-structure, 
i.e., a coframing~$(\theta_0,\eta_0)$ satisfying
\be\label{eq: initstrcteqns}
\begin{aligned}
\d\theta_0 &= \iC\,\eta_0\w\overline{\eta_0}\,,\\
\d\eta_0 &= 2\iC\,\theta_0\w\bigl(a_0\,\eta_0 + b_0\,\overline{\eta_0}\bigr)
\end{aligned}
\ee
for some constants~$a_0 = \overline{a_0}$ and~$b_0$.  

Of course, this coframing is real-analytic in appropriate
coordinates and hence the realization of~$M^3$ as a nondegenerate
hypersurface in a unimodular complex surface is locally unique.
In this particular case, one can explicitly define the complex
structure on a neighborhood of~$M\times\{0\}\subset M\times\bbR$
that realizes this with a holomorphic volume form 
\be
\Upsilon = (\theta + \iC\,\d t) \w \eta
\ee
where~$\theta$ and~$\eta$ are $t$-dependent $1$-forms on~$M$.

Computation yields that, when regarded as $1$-forms
on an open set in~$M\times\bbR$, the $1$-forms~$\theta$ and~$\eta$ 
must satisfy
\be\label{eq: thetaetatdep streqs}
\begin{aligned}
\d\theta &= \iC\,\eta\w\overline{\eta} + {\ts\frac43}\,a\,\theta\w\d t\,,\\
\d\eta &= 2\iC\,(\theta+\iC\,\d t)\w\bigl(a\,\eta + b\,\overline{\eta}\bigr)
              - {\ts\frac43}\,a\,\eta\w\d t
\end{aligned}
\ee
where~$a=\overline a$ and~$b$ are functions of~$t$ that satisfy
the differential equations
\be\label{eq: ab ode}
\begin{aligned}
a'(t) &= {\ts\frac43}\,a(t)^2 + 4\,b(t)\,\overline{b(t)}\\
b'(t) &= {\ts\frac{16}3}\,a(t) b(t)
\end{aligned}
\ee
and initial conditions~$a(0) = a_0$ and~$b(0) = b_0$. 

Note that~\eqref{eq: ab ode} does not have any periodic 
solutions other than the fixed point~$(a,b)\equiv (0,0)$
since, except in this case, $a'$ is strictly positive.

Thus, if~$I_0\subset\bbR$ is the maximal interval containing~$0\in\bbR$
on which~$a$ and~$b$ can be defined satisfying~\eqref{eq: ab ode}
and the given initial conditions, then Cartan's standard existence
theorem implies that, for any $3$-dimensional~$M_0$ 
endowed with a coframing~$(\theta_0,\eta_0)$ satisfying
\eqref{eq: initstrcteqns}, there will exist a time-dependent
coframing~$(\theta,\eta)$ on~$M_0\times I_0$ 
satisfying~\eqref{eq: thetaetatdep streqs} where~$a$ and~$b$
satisfy~\eqref{eq: ab ode} and where~$(\theta,\eta)$ pull back
to~$M_0\times\{0\}$ to become~$(\theta_0,\eta_0)$.

In fact, though, one does not need to invoke Cartan's theorem, 
since one can explicitly find~$(\theta,\eta)$ in the
form
\be
\theta = \bigl(\det F(t)\bigr)\,\theta_0
\qquad\text{and}\qquad
\bpm\eta\\ \overline{\eta}\epm
 = F(t)\,\bpm\eta_0\\ \overline{\eta_0}\epm
\ee
where~$F(t)$, defined for~$t\in I_0$, 
satisfies the initial condition~$F(0)=I_2$ and
the differential equation
\be
F'(t) = -2\bpm \frac13 a(t) & b(t)\\[3pt]
           \overline{b(t)} & \frac13 a(t)\epm
\,F(t).
\ee
where~$\bigl(a(t),b(t)\bigr)$ is the solution on~$I_0$ of~\eqref{eq: ab ode}
that satisfies~$a(0) = a_0$ and~$b(0) = b_0$.

Now, the general solution of \eqref{eq: ab ode} is
expressed in elliptic functions.  However, the problem can
be simplified and a partial explicit integration effected:
For any solution of~\eqref{eq: ab ode}, either~$b(t)$ vanishes
identically or it never vanishes.  

When~$b(t)\equiv0$, the remaining equation for~$a$ is integrated as
\be\label{eq: atwhenb=0}
a(t) = \frac{a_0}{1-\frac43 a_0 t}.
\ee
It then follows that
\be
\theta = \bigl(1-{\ts\frac43}a_0 t\bigr)\,\theta_0
\qquad\text{and}\qquad
\eta = \bigl(1-{\ts\frac43}a_0 t\bigr)^{\frac12}\,\eta_0\,.
\ee
When~$a_0>0$, this gives the evolution of
the sphere of radius~${a_0}^{-3/4}$
in~$\bbC^2$ (endowed with its standard volume form).
 
\begin{example}[Circle bundles over~$K=-1$ surfaces]
\label{ex: homa<0b=0evol} 
This is a continuation of the discusion begun 
in Example~\ref{ex: homa<0b=0}.
Recall that the hypersurfaces~$M_r\subset T^*\bbD$ 
constructed there have~$a = -\frac12 r^{-2/3}$.  It follows 
from~\eqref{eq: atwhenb=0}
that the unimodular normal flow with~$M_{r_0}$ as
initial hypersurface is described by~$M_{r(t)}$ where
\be
r(t) = (r_0^{2/3}+{\ts\frac23}\,t)^{3/2}.
\ee
Of course, this flow exists for all positive time, but 
exists only for a finite time in the past.  Note that
as~$t$ decreases towards the singularity at~$t=-\frac32\,r_0^{2/3}$,
the hypersurface~$M_{r(t)}$ converges to the zero section of~$X\to\bbD$,
i.e., the hypersurface collapses onto a complex curve.  
\end{example}

Assume, henceforth, that~$b(t)$ is nowhere-vanishing.  Then one 
finds that the ratio~$\overline{b}/b$ is constant and that,
after applying scaling as discussed in Remark~\ref{rem: scaleffects}, 
all solutions can be deduced from solutions in which~$b(t)$ is positive
and real, so assume this from now on.  Then the obvious homogeneity
of~\eqref{eq: ab ode} implies that there must be a homogeneous
first integral and, indeed, one finds that the function $(a^2-b^2)^2/b$
is constant.  In other words, each integral curve of~\eqref{eq: ab ode}
with~$b$ positive and real lies in a quartic plane curve of the form
\be\label{eq: abquarticcurve}
(a^2-b^2)^2 - \lambda^3\, b = 0
\ee
for some constant~$\lambda\ge0$.

When~$\lambda = 0$, one has the solutions
\be
a(t) = \frac{a_0}{1-\frac{16}3 a_0 t}
\qquad\text{and}\qquad
b(t) = \frac{b_0}{1-\frac{16}3 a_0 t}
\ee
with~$a_0=\pm b_0\not=0$.

When~$\lambda>0$, the curve~\eqref{eq: abquarticcurve} is 
irreducible over~$\bbC$ and of genus~$1$.  It is smooth
in the finite part of the plane and has two ordinary double
points on the line at infinity.  The functions~$a$ and $b$
each have four simple poles on the normalized curve, representing
the four points on the line at infinity.  The function~$b$ has
a quadruple zero at~$(a,b) = (0,0)$ and no other zeros, while
the function~$a$ has four simple zeros, only two of which, those
at~$(a,b) = (0,0)$ and~$(a,b) = (0,\lambda)$, are real.  
The curve~\eqref{eq: abquarticcurve} has two real branches:  One
that passes through~$(a,b) = (0,0)$ and that, except for this point,
lies in the union of the sectors defined by~$0<b<-a$ and~$0<b<a$
and one that passes through~$(a,b) = (0,\lambda)$ and lies
in the sector~$b>|a|$.
	Away from the point~$(a,b)=0$, one has
\be
\d t = \frac{3\,\d b}{16 a b}
     = \frac{3\,\d a}{4(a^2 + 3 b^2)},
\ee  
implying that~$\d t$ is a meromorphic differential with
a double pole at~$(a,b)=0$ and no other zeros or poles on the
real branches (note that~$\d b$ vanishes at~$(a,b)=(0,\lambda)$),
even on the line at infinity.

	In particular, the integral of $\d t$ over the branch 
in the sector~$b>|a|$ is finite, so that these solutions exist
for only a finite time both forward and backwards.
	
	Similarly, any solution of~\eqref{eq: abquarticcurve}
that lies in the sector~$0<b<-a$ exists for only a finite time
in the past but an infinite time in the future, while any solution
of~\eqref{eq: abquarticcurve} that lies in the sector~$0<b<a$ 
exists for an infinite time in the past but only a finite time
in the future.

	Note that, since~$\bigl(\log\bigl(\det F(t)\bigr)\bigr)'
=-{\ts\frac43}a(t) = -{\ts\frac14}\,\bigl(\log b(t)\bigr)'$,
it follows that
\be
\det F(t) = \bigl(b(t)/b_0\bigr)^{-\frac14},
\ee
so
\be
\theta = \bigl(b(t)/b_0\bigr)^{-\frac14}\,\theta_0\,.
\ee
\end{example}

\begin{example}[Explicit homogeneous hypersurface flows]
Finally, consider the flow for the homogenous nondegenerate
hypersurfaces~$M_{h,q}\subset X_q$ in either Example~\ref{ex: SO3exs}
(the $\SO(3)$-homogeneous examples) or Example~\ref{ex: PSL2exs}
(the $\PSL(2,\bbR)$-homogeneous examples).

In either case, the evolution equations~\eqref{eq: ab ode}
for the invariants~$(a,b)$ coupled with the 
equations~\eqref{eq: abasrxi} or~\eqref{eq: abashq}
imply that the hypersurface~$M_{h_0,q}$ will evolve
as the family~$M_{h(t),q}$ where~$h(t)$ satisfies the
ordinary differential equation initial value problem
\be\label{eq: h ode}
h'(t) = 2\bigl(|q|^2-h(t)^2\bigr)^{1/3},
\qquad
h(0) = h_0\,.
\ee

If~$h_0 > |q|$, which happens in the $\SO(3)$-invariant
examples (and only in those examples), it follows that 
the flow exists for an infinite time into the
past (and~$h(t)\to\infty$ as~$t\to-\infty$) 
but only a finite time into the future, 
as $h(t)$ will go down to~$|q|$ in time
\be
T(h_0,q) = \int_{|q|}^{h_0} \frac{\d h}{2\bigl(h^2-|q|^2\bigr)^{1/3}}
<\infty.
\ee
In particular~$M_{h(t),q}$ will converge in finite time to~$M_{|q|,q}
\subset X_q$, which is a point if~$q=0$%
\footnote{The reader may wonder what would happen if one were
to replace the singular unimodular quadric~$(X_0,\Upsilon_0)$ by its
crepant resolution~$(\hat X_0,\hat\Upsilon_0)$.  The answer, not
surprisingly, is that~$M_{h(t),0}$, regarded as a hypersurface 
in~$\hat X_0$, would collapse in finite time 
to the exceptional curve in~$\hat X_0$.}
and an embedded, totally real $2$-sphere in~$X_q$ if~$|q|>0$.

If~$h_0 < -|q|$, which happens in the $\PSL(2,\bbR)$-invariant
examples (and only in those examples), it follows that
the flow exists for all future time 
(and~$h(t)\to-\infty$ as~$t\to\infty$)  
but extends only a finite time into the past, 
as $h(t)$ will have come down from a starting
value of~$-|q|$ at the (past, i.e., negative) time
\be
T(h_0,q) 
= \int^{-|q|}_{h_0} \frac{\d h}{2\bigl(|q|^2-h^2\bigr)^{1/3}}
> -\infty.
\ee

Recall that~$M_{h,q}$ is the disjoint union of~$M^+_{h,q}$
and~$M^-_{h,q}$ (as defined in Example~\ref{ex: PSL2exs}). 

It is interesting to consider what happens to each of these
smooth hypersurfaces as~$h$ approaches~$|q|$ from below.
When~$q=0$, the hypersurfaces~$M^\pm_{h,0}$ as~$h\to 0^-$
converge to~$M_{0,0}$ in the sense that every point of~$M_{0,0}$
is a limit of a sequence~$z_k\in M^+_{h_k,0}$ with~$h_k\to 0^-$,
with a similar statement for the hypersurfaces~$M^-_{h,0}$.
This is so in spite of the fact that, under the actual unimodular
normal flow, all of the points of~$M_{h,0}$ with~$h<0$ flow to
the origin itself.  Thus, there is a significant difference
between the limit of the pointwise flow and the `hypersurface
flow'.  When~$q\not=0$, the situation is even more interesting.
As $h\to -|q|^-$, the hypersurface~$M^+_{h,q}$ converges to 
a proper subset of~$M_{-|q|,q}$ that consists of the singular
locus~$\iC\,H_{-1}$ (a surface) and two of the four components
of the smooth locus of~$M_{-|q|,q}$, say,~$M^{+,+}_{-|q|,q}$ 
and~$M^{-,-}_{-|q|,q}$.  In particular, the union of these three
$\PSL(2,\bbR)$-orbits is a singular, degenerate hypersurface, 
with singular locus equal to the $2$-dimensional orbit~$\iC\,H_{-1}$.
On the other hand, as $h\to -|q|^-$, the hypersurface~$M^-_{h,q}$ 
converges to a proper subset of~$M_{-|q|,q}$ that consists 
of the singular locus~$\iC\,H_{-1}$ (a surface) and two of 
the four components of the smooth locus of~$M_{-|q|,q}$, 
say,~$M^{+,-}_{-|q|,q}$ and~$M^{-,+}_{-|q|,q}$.  In particular, 
the union of these three $\PSL(2,\bbR)$-orbits is a singular, 
degenerate hypersurface, with singular locus equal to the 
$2$-dimensional orbit~$\iC\,H_{-1}$.  In each case, under the
actual unimodular normal flow, the points of~$M_{h,q}$ with~$h<-|q|$
converge to the $2$-dimensional singular locus~$\iC\,H_{-1}$, 
while the hypersurface itself converges to a (singular) 
$3$-dimensional hypersurface.

Finally, consider what happens to the 
hypersurfaces~$M_{h_0,q}\subset X_q$ with~$|h_0|<|q|$. 
Again, these are $\PSL(2,\bbR)$-invariant hypersurfaces.
The unimodular normal flow with this initial hypersurface
$M_{h(t),q}$ exists for~$T^-(h_0,q)<t<T^+(h_0,q)$, 
where
\be
T^+(h_0,q) 
= \int^{|q|}_{h_0} \frac{\d h}{2\bigl(|q|^2-h^2\bigr)^{1/3}} 
< +\infty
\ee
and
\be
T^-(h_0,q)
= \int^{-|q|}_{h_0} \frac{\d h}{2\bigl(|q|^2-h^2\bigr)^{1/3}}
>-\infty.
\ee

Several things are interesting about these flows:  
First, the hypersurface~$M_{0,q}$ is CR-flat (since these have~$ab=0$)
but~$M_{h,q}$ is not CR-flat when~$hq\not=0$.  Thus, this furnishes
an example of a CR-flat hypersurface for which 
the unimodular flow does not preserve CR-flatness.  
Second, as $h$ increases to the value~$|q|$, the hypersurfaces~$M^+_{h,q}$
converge to the totally real surface~$M^+_{|q|,q}\subset X_q$,
i.e., the upper nappe of the (real) unit hyperboloid of two sheets
while the hypersurfaces~$M^-_{h,q}$ converge to the totally real 
surface~$M^-_{|q|,q}\subset X_q$, i.e., the lower nappe of the (real) 
unit hyperboloid of two sheets.  Thus, the `singularity' that develops
in forward time from~$M_{h,q}$ under the unimodular normal flow 
is collapse onto a totally real surface.
Third, as $h$ decreases to the value~$-|q|$, the hypersurfaces~$M^+_{h,q}$
converge to the singular real hypersurface in~$X_q$ that
is the union of the three $\PSL(2,\bbR)$-orbits~$M^{+,+}_{-|q|,q}$,
$M^{-,+}_{-|q|,q}$, and~$\iC\,H_{-1}$ while the hypersurfaces~$M^-_{h,q}$ 
converge to the singular real hypersurface in~$X_q$ that
is the union of the three $\PSL(2,\bbR)$-orbits~$M^{+,-}_{-|q|,q}$,
$M^{-,-}_{-|q|,q}$, and~$\iC\,H_{-1}$.  (Nevertheless, under this
flow, each point of $M_{h,q}$ converges to a point of~$\iC\,H_{-1}$.)
\end{example}

\section{Invariants of $3$-dimensional CR-manifolds}
\label{sec: CRinvs}

This section contains an exposition 
of Cartan's solution~\cite{ECartanCR1, ECartanCR2} in 1932
of the equivalence problem for nondegenerate hypersurfaces 
in complex $2$-manifolds (with no volume specified).
The main reason for including it here is that it will be
used to compute the invariants of the underlying CR-structure
associated to the canonical coframing of a pseudoconvex
real hypersurface in a unimodular complex surface.

This equivalence problem is much more subtle than 
the case of hypersurfaces in a unimodular surface.
For another exposition of Cartan's solution in more `modern' 
language, the reader might compare 
Jacobowitz'~\cite[Chapters 5-7]{MR93h:32023}, though the
notation is different from that of this article.

Of course, the theory has been extensively developed since Cartan's work, 
with the general solution for a nondegenerate hypersurface in a complex 
$n$-manifold being the subject of a famous paper 
by Chern and Moser~\cite{MR5413112},
and an earlier paper by Tanaka~\cite{MR263086}.

\subsection{The geometric problem and its $G$-structure}
\label{ssec: geomprobGstr}
Suppose that~$M^3\subset X$ is a smooth real hypersurface 
in a complex $2$-manifold~$X$, which, since all the 
considerations are local, can be taken to be~$\bbC^2$ if desired.

\subsubsection{The notion of a \textrm{CR}-structure} 
For each~$x\in M$, the tangent plane~$T_xM$ 
cannot be a complex subspace of $T_xX$, but contains 
a unique complex subspace~$D_x\subset T_xM$ of complex dimension~$1$.  
Thus, $M$ inherits a geometric structure from being immersed 
as a hypersurface in a complex $2$-manifold.

\begin{definition}
A (smooth) \emph{\textrm{CR}-structure} on a $3$-manifold~$M$ is a choice of a 
(smooth) rank~$2$ subbundle~$D\subset TM$ together with a choice of complex 
structure on~$D$, i.e., a smooth bundle map~$J:D\to D$ satisfying~$J^2 = 
-\text{Id}_D$.
\end{definition}

In the real-analytic category, every CR-structure on a $3$-manifold 
is locally induced by an immersion into~$\bbC^2$.

\begin{proposition}[Local realization of analytic CR-structures]
Let $\bigl(D,J\bigr)$ be a real-analytic \textrm{CR}-structure on~$M^3$.  
Then for each point~$x\in M$ there exists an $x$-neighborhood~$U$ 
and a real-analytic embedding~$Z: U\to\bbC^2$ so that $\bigl(D,J\bigr)$ 
is the \textrm{CR}-structure on~$U$ induced by the embedding~$Z$.
\end{proposition}

\begin{proof} 
On a neighborhood~$U$ of~$x$ choose a real-analytic, nonvanishing 
real $1$-form~$\rho$ that annihilates~$D$ and a real-analytic, complex valued 
$1$-form~$\eta$ linearly independent from~$\rho$ that satisfies~$\eta(Jv) = 
\iC\,\eta(v)$ for all~$v\in D$.  Then any complex-valued $1$-form~$\zeta$ 
on~$U$ that satisfies~$\zeta(Jv) = \iC\,\zeta(v)$ is a linear combination 
of~$\rho$ and~$\eta$.  

To construct the desired~$Z$, it suffices to find two complex 
functions~$z^1$ and $z^2$ in a neighborhood of~$x$ 
whose differentials are linearly independent and that satisfy~$\d z^k(Jv) = 
\iC\,\d z^k(v)$, i.e., so that $\d z^k\w\rho\w\omega = 0$.  

Now, on~$N=U\times\bbC$ with second projection~$z:N\to\bbC$, 
let~$\cI$ be the ideal generated by the two $3$-forms 
that are the real and imaginary parts of~$\d z\w\rho\w\omega$.  
The characters of~$\cI$ are~$s_i=0$ for $i\not=2$ and $s_2=2$.  
Meanwhile, the space of $3$-dimensional integral elements of~$\cI$ 
that satisfy the independence condition~$\rho\w\omega\w\bar\omega\not=0$ 
is of dimension~$4$. Thus, the system~$\cI$ is in involution. 

Choose two integral manifolds~$\Sigma_i$, $i=1,2$ 
of~$\cI$ that pass through~$(x,0)\in N$ but that are not tangent there.  
Each is then the graph of a function~$z^i$ 
that satisfies~$\d z^k\w\rho\w\omega = 0$ 
and the condition that the two integral manifolds 
not be tangent is equivalent to $\d z^1\w \d z^2\not=0$.
\end{proof}

\begin{remark}[The need for the analyticity hypothesis]
The famous Levi-Nirenberg example~\cite{MR83f:53022, MR56:9048} 
shows that the assumption of real-analyticity is necessary here.
\end{remark}

\subsubsection{A $G$-structure associated to a \textrm{CR}-structure} 
Suppose now that~$M^3$ is endowed with a $CR$-structure~$\bigl(D,J\bigr)$.  
Let~$V=\bbR\oplus\bbC$ and think of~$V$ as the space of columns of height~$2$ 
whose first entry is real and whose second entry is complex.  

A coframe~$u: T_xM\to V$ will be said to be \emph{$0$-adapted} 
to~$\bigl(D,J\bigr)$ if $u(D_x) = \bbC\subset V$ and, 
moreover, $u(Jv) = i\,u(v)$ for all~$v\in D_x$.  
Let~$B_0\subset F^*(M,V)$ denote the space of~$0$-adapted $V$-valued 
coframes on~$M$.  

If $u$ and $u^*$ lie in~$B_0$ and share the same basepoint, then
\be
u^* = \begin{pmatrix} r&0\\ b&a\\ \epm\>u,
\ee
where $r$ is a real number and $a$ and $b$ are complex, with $a\not=0$.  
Thus, $B_0$ is a $G_0$-structure on~$M$ where
\be
G_0 = \left\{\ \bpm r&0\\ b&a\\ \epm
\ \vrule\ r\in\bbR^*,\ a\in\bbC^*,\ \text{and}\ b\in\bbC\ \right\}.
\ee

Conversely, given a $G_0$-structure~$B_0$ on~$M$, there is canonically 
associated to it a unique CR-structure~$\bigl(D,J\bigr)$ that gives rise 
to it via this construction.  
Thus, the two sorts of structures are equivalent.

\subsection{The first analysis}\label{ssec: 1stanalysis}
Now let~$B_0$ be a $G_0$-structure on~$M^3$.  Write the 
canonical~$V$-valued $1$-form~$\omega$ on~$B_0$ in the form
\be
\omega = \bpm \theta\\ \eta\\ \epm
\ee
where~$\theta$ is a real-valued $1$-form 
and $\eta$ is a complex-valued $1$-form.  

\subsubsection{The first structure equation}
The first structure equation can be written in the form
\be
\d\bpm \theta\\ \eta\\ \epm
=-\bpm \rho_0&0\\ \beta_0 & \alpha_0\\ \epm
\w\bpm \theta\\ \eta\\ \epm
+\bpm \theta\w\bigl(b\,\eta+\bar b\,\bar\eta\bigr) + \iC 
L\,\eta\w\bar\eta\\
\theta\w\bigl(c\,\eta+ e\,\bar\eta\bigr) + T\,\eta\w\bar\eta\\ \epm
\ee
where $L$ is a real function on~$B_0$ but the other coefficients 
are allowed to be complex.  Clearly, by adding multiples of~$\theta$, 
$\eta$ and $\bar\eta$ to the pseudo-connection forms~$\rho_0$, $\alpha_0$, 
and $\beta_0$, it can be arranged that $b=c=e=T=0$, but~$L$ cannot
be affected by such changes.  

Thus, the structure equations can be assumed to have the form
\be
\d\bpm \theta\\ \eta\\ \epm
=-\bpm \rho_0&0\\ \beta_0 & \alpha_0\\ \epm
\w\bpm \theta\\ \eta\\ \epm
+\bpm \iC L\,\eta\w\bar\eta\\0\\ \epm.
\ee
Differentiating the first equation~$\d\theta=-\rho_0\w\theta 
+ \iC L\,\eta\w\bar\eta$ and reducing modulo~$\theta$ gives the relation
\be
\d L\equiv L\bigl(\alpha_0+\bar\alpha_0-\rho_0\bigr)\mod \theta,\eta,\bar\eta,
\ee
so, on a given fiber of~$B_0$, 
either~$L$ vanishes identically or is nowhere zero there. 

\subsubsection{The case~$L\equiv0$}
The case where~$L$ vanishes identically, i.e., the intrinsic torsion of the 
$G_0$-structure vanishes, turns out not to be very interesting.  In this case, 
one can calculate that the characters of the Lie algebra~$\eug_0$ are $s_1=3$, 
$s_2=1$, and $s_3 =0$.  Moreover the variability of the pseudo-connection is of 
dimension~$5=s_1+2s_2$, so $G_0$ is semi-involutive and all of the real-analytic 
$G_0$-structures with vanishing torsion are equivalent.  Thus, it makes sense to 
concentrate on the (generic) case where $L$ is nowhere-vanishing.%
\footnote{The `intermediate' case, in which~$L$ vanishes on 
some proper closed subset has also been intensively studied,
this body of work will play no role in this article.}

\subsubsection{The case~$L\not=0$}
Now, there is a direct geometric interpretation of~$L$.  
Since~$\theta$ is a nonzero multiple of~$\pi^*(\sigma)$,
where $\sigma$ is any nonvanishing $1$-form with $D=\ker\sigma$, 
it follows that $\theta\w \d\theta = \iC L\,\theta\w\eta\w\bar\eta$ 
is nonzero if and only if $\sigma\w \d\sigma$ is nonzero, 
i.e., if and only if~$D$ is a contact plane field on~$M^3$.

\begin{definition} 
A CR-structure~$\bigl(D,J\bigr)$ on~$M^3$ is \emph{nondegenerate} 
if $D$ is nowhere-integrable, i.e., is a contact structure on~$M$.
\end{definition}

Thus, the condition that~$L$ be nowhere-vanishing is the condition that the 
CR-structure be nondegenerate.  From now on, this will be assumed 
to be the case.  

\subsubsection{The first structure reduction}
This assumption leads directly to the first reduction:  Set
\be
B_1 = \{\ u\in B_0\,\vrule\, L(u) = 1\ \}.
\ee
Then~$B_1$ is a $G_1$-structure on~$M$ where
\be
G_1 = \left\{\ \bpm a\bar a&0\\ b&a\\ \epm
\ \vrule\ a\in\bbC^*\ \text{and}\ b\in\bbC\ \right\}.
\ee
Pulling all of the forms on~$B_0$ back to~$B_1$ 
and giving them the same names, 
the structure equations on~$B_1$ now read
\be
\d\bpm \theta\\ \eta\\ \epm
=-\bpm \alpha_0+\bar\alpha_0&0\\ \beta_0 & \alpha_0\\ \epm
\w\bpm \theta\\ \eta\\ \epm
+\bpm (\alpha_0+\bar\alpha_0-\rho_0)\w\theta+\iC\,\eta\w\bar\eta\\0\\ 
\epm.
\ee
where~$\alpha_0+\bar\alpha_0-\rho_0 = a\theta + b\eta + \bar b\bar\eta$ 
for some functions~$a$ and $b$ on~$B_1$.  
Subtracting~$b\eta$ from~$\alpha_0$ reduces the function~$b$ 
to zero and the structure equations become
\be\label{eq: 1streqsonB_1}
\d\bpm \theta\\ \eta\\ \epm
=-\bpm \alpha_0+\bar\alpha_0&0\\ \beta_0 & \alpha_0\\ \epm
\w\bpm \theta\\ \eta\\ \epm
+\bpm \iC\,\eta\w\bar\eta\\0\\ \epm.
\ee
Now the torsion is constant.  

If the algebra~$\eug_1$ were involutive, then reaching this point 
would imply that any two nondegenerate $G_0$-structures were locally 
equivalent.  However, one easily computes that the characters of 
this algebra are~$s_1=3$, $s_2=1$, and $s_3 = 0$ while the 
pseudo-connections with this torsion are determined up to a replacement 
of the form~$\bigl(\alpha_0,\beta_0\bigr)
\mapsto \bigl(\alpha_0^*,\beta_0^*\bigr)$ where
\be
\bpm \alpha_0^*\\ \beta_0^*\\ \epm
=\bpm \alpha_0\\ \beta_0\\ \epm
+\bpm s^1&0\\ s^2&s^1\\ \epm
\bpm \theta \\ \eta \\ \epm,
\ee
and $s^1$ and $s^2$ are arbitrary complex-valued functions on~$B_1$.
Thus, $\dim\eug_1^{(1)}=4<s_1+2s_2+3s_3=5$, so~$\eug_1$ is not involutive.  
Hence, there remains the possibility that there will be differential 
invariants at some higher order.

\subsection{Prolongation and further reductions}\label{ssec: prolong&red}
According to the prescription of the method of equivalence, 
I now construct a $\eug_1^{(1)}$-bundle~$B^{(1)}_1$ over $B_1$ 
that consists of the coframes on~$B_1$ with values in $V\oplus\eug_1$ 
that satisfy the structure equations of~$B_1$.  

For simplicity, I will identify $V\oplus\eug_1$ 
with~$\bbR\oplus\bbC^3$, thought of as the columns of height~$4$ 
with the first entry real and the remaining three complex.  

In the trivialization $B^{(1)}_1 = B_1\times\eug^{(1)}_1$ 
induced by the section~$B_1\to B^{(1)}_1$ 
represented by a choice of~$\alpha_0$ and~$\beta_0$ on~$B_1$ 
satisfying~\eqref{eq: 1streqsonB_1}, 
the canonical $1$-form~$\omega^{(1)}$ has the form
\be
\omega^{(1)} = \bpm\theta\\ \eta\\ \alpha\\ \beta\\ \epm
=\bpm1&0&0&0\\ 0&1&0&0\\ s^1&0&0&0\\ s^2&s^1&0&0\\ \epm^{-1}
\bpm\theta\\ \eta\\ \alpha_0\\ \beta_0\\ \epm
=\bpm\theta\\ \eta\\ \alpha_0-s^1\theta\\ \beta_0-s^2\theta-s^1\eta\\
\epm
\ee
where, of course, the functions~$s^1$ and $s^2$ now represent coordinates 
on~$\eug^{(1)}_1$ and so are independent from the functions on~$B_1$.  
The structure equations on~$B^{(1)}_1$ have the form:
\be\label{eq: 1streqsB1_1}
\d\bpm\theta\\ \eta\\ \alpha\\ \beta\\ \epm
=-\bpm0&0&0&0\\ 0&0&0&0\\ \sigma_0^1&0&0&0\\ 
\sigma_0^2&\sigma_0^1&0&0\\ \epm
\bpm\theta\\ \eta\\ \alpha\\ \beta\\ \epm
+\bpm
-(\alpha+\bar\alpha)\w\theta + \iC\,\eta\w\bar\eta\\
-\beta\w\theta-\alpha\w\eta\\
T_\alpha\\ T_\beta\\ \epm
\ee
where~$T_\alpha$ and $T_\beta$ represent the torsion terms 
associated to those components of the canonical $1$-form 
while $\sigma^1_0$ and $\sigma^2_0$ are $1$-forms 
that satisfy~$\sigma^i_0\equiv \d s^i$ modulo forms semi-basic 
for the projection~$B^{(1)}_1\to B_1$ but that are otherwise arbitrary.

Computing the exterior derivatives of the first two equations 
of~\eqref{eq: 1streqsB1_1} yields
\be
\begin{aligned}
0 & {} = \d\bigl(\d\theta\bigr) = {} %&
-\bigl(T_\alpha +\overline{T_\alpha} - \iC\,\beta\w\bar\eta
+\iC\,\bar\beta\w\eta\bigr)\w\theta,\\
0 & {} = \d\bigl(\d\eta\bigr)   = {} %&
-\bigl(T_\beta +\beta\w\bar\alpha\bigr)\w\theta
-\bigl(T_\alpha+\iC\,\beta\w\bar\eta\bigr)\w\eta,\\
\end{aligned}
\ee

Setting $T_\alpha^* = T_\alpha + \iC\,\beta\w\bar\eta 
+ 2\iC\,\bar\beta\w\eta$ and $T_\beta^* = T_\beta +\beta\w\bar\alpha$, 
these equations can be written in the form
\be
\bigl(\,T^*_\alpha +\overline{T^*_\alpha}\,\bigr)\w\theta
= \bigl(T^*_\beta\bigr)\w\theta+\bigl(T^*_\alpha\bigr)\w\eta = 0
\ee
and the second of these equations implies, via Cartan's Lemma, 
that there exist $1$-forms~$\psi_1$, $\psi_2$, and $\psi_3$ so that
\be
\begin{aligned}
T^*_\alpha &= \psi_2\w\theta +\psi_1\w\eta\\
T^*_\beta  &= \psi_3\w\theta +\psi_2\w\eta\\
\end{aligned}
\ee
Since~$T^*_\alpha$ and $T^*_\beta$ are semi-basic, the $\psi_i$ 
must be also.  Thus, by subtracting~$\psi_2$ from~$\sigma^1_0$ and $\psi_3$ 
from~$\sigma^2_0$, it can be arranged that~$\psi_2=\psi_3=0$.  
Then the remaining equation on~$T^*_\alpha$ becomes
\be
\bigl(\,\psi_1\w\eta + \overline{\psi_1\w\eta}\,\bigr)\w\theta = 0,
\ee
which implies that~$\psi_1\w\eta = b\,\eta\w\theta + R\,\eta\w\bar\eta$, 
where~$b$ is a complex function and~$R$ is a real function.  
By adding $b\,\eta$ to~$\sigma^1_0$, it can be arranged that~$b=0$, 
so that the structure equations now take the form
\be
\d\bpm\theta\\ \eta\\ \alpha\\ \beta\\ \epm
=-\bpm0&0&0&0\\ 0&0&0&0\\ \sigma_0^1&0&0&0\\ 
\sigma_0^2&\sigma_0^1&0&0\\ \epm
\bpm\theta\\ \eta\\ \alpha\\ \beta\\ \epm
+\bpm
-(\alpha+\bar\alpha)\w\theta + \iC\,\eta\w\bar\eta\\
-\beta\w\theta-\alpha\w\eta\\
- \iC\,\beta\w\bar\eta - 2\iC\,\bar\beta\w\eta
+R\,\eta\w\bar\eta \\ -\beta\w\bar\alpha\\ \epm.
\ee

Now, computing the exterior derivative of the $\d\alpha$ 
equation modulo~$\theta$ yields
\be
0 = \d\bigl(\d\alpha\bigr)
\equiv \bigl(dR - (\alpha+\bar\alpha)R 
-2\iC(\overline{\sigma^1_0}-\sigma^1_0)\ \bigr)\w\eta\w\bar\eta
\mod \theta,
\ee
which implies
\be\label{eq: dRcong}
\d R \equiv (\alpha+\bar\alpha)R 
+2\iC(\overline{\sigma^1_0}-\sigma^1_0)
\mod \theta,\eta,\bar\eta.
\ee
In particular, on each fiber of $B^{(1)}_1\to B_1$, the relation~$\d R = 
2\iC\,d(\overline{s^1}-s^1)$ holds.  

It follows that the equation~$R=0$ 
defines a $G_2$-structure~$B_2\subset B^{(1)}_1$ on~$B_1$, 
where $G_2$ is the subgroup consisting of those matrices 
in~$\eug^{(1)}_1$ for which~$s^1$ is real. 
 
Now pull back all of the forms and functions on~$B^{(1)}_1$ to~$B_2$, 
write~$\sigma^1_0 = \sigma_0 + \iC\,\tau$ where $\sigma_0$ and $\tau$ 
are real $1$-forms, and write the structure equations on~$B_2$ in the form
\be
\d\bpm\theta\\ \eta\\ \alpha\\ \beta\\ \epm
=-\bpm0&0&0&0\\ 0&0&0&0\\ \sigma_0&0&0&0\\ 
\sigma_0^2&\sigma_0&0&0\\ \epm
\bpm\theta\\ \eta\\ \alpha\\ \beta\\ \epm
+\bpm
-(\alpha+\bar\alpha)\w\theta + \iC\,\eta\w\bar\eta\\
-\beta\w\theta-\alpha\w\eta\\
- \iC\,\beta\w\bar\eta - 2\iC\,\bar\beta\w\eta
-\iC\,\tau\w\theta \\ 
-\beta\w\bar\alpha-\iC\,\tau\w\eta\\ \epm.
\ee
The congruence~\eqref{eq: dRcong} now implies that~$\tau 
= a\,\theta + b\,\eta + \bar b\,\bar\eta$ for some real 
function~$a$ and complex function~$b$ on~$B_2$.
By adding~$\iC a\,\eta$ to~$\sigma^2_0$, it can be arranged that~$a=0$, 
but~$b$ cannot be absorbed.  

The last two structure equations now read
\be
\begin{aligned}
\d\alpha &= -\sigma_0\w\theta - \iC\,\beta\w\bar\eta - 
2\iC\,\bar\beta\w\eta-\iC\,(b\,\eta + \bar b\,\bar\eta)\w\theta,\\
\d\beta &= -\sigma^2_0\w\theta-\sigma_0\w\eta 
-\beta\w\bar\alpha+\iC\,\bar b\,\eta\w\bar\eta,\\
\end{aligned}
\ee
and it remains to determine how~$b$ varies on the fibers of~$B_2\to B_1$.  

To do this, compute the exterior derivative 
of the first of these equations and write it in the form
\be\label{eq: ddalpha}
\begin{aligned}
0 = \d(\d\alpha) = {} 
&-\left(\d\sigma_0-(\alpha+\bar\alpha)\w\sigma_0 -\iC\,\beta\w\bar\beta
+{\ts{\frac12}}\iC
\bigl(\sigma^2_0\w\bar\eta-\overline{\sigma^2_0}\w\eta\bigr)\,\right)\w\theta\\
&-\iC\,\left(db -(2\alpha+\bar\alpha)\,b 
+{\ts{\frac32}}\,\overline{\sigma^2_0}\,\right)\w\eta\w\theta\\
&-\iC\,\left(\d\bar b -(2\bar\alpha+\alpha)\,\bar b 
+{\ts{\frac32}}\,\sigma^2_0\,\right)\w\bar\eta\w\theta\\
\end{aligned}
\ee
The imaginary part of~\eqref{eq: ddalpha} implies that 
\be
\d b \equiv (2\alpha+\bar\alpha)\,b -{\ts{\frac32}}\,\overline{\sigma^2_0}
\mod \theta,\eta,\bar\eta,
\ee
which implies that, on each fiber of~$B_2\to B_1$, 
an equation of the form~$\d b = {\frac32}\,\overline{\d s^2}$ holds.  

In particular, the equation~$b=0$ defines a~$G_3$-structure~$B_3\subset B_2$ on~$B_1$, where~$G_3$ is the 1-dimensional subgroup of~$G_2$ 
defined by the equation~$s^2=0$.  

Now pull back all of the forms and functions involved to~$B_3$.  
The structure equations take the form
\be
\d\bpm\theta\\ \eta\\ \alpha\\ \beta \epm
=
-\bpm0&0&0&0\\ 0&0&0&0\\ \sigma_0&0&0&0\\ 0&\sigma_0&0&0 \epm
\bpm\theta\\ \eta\\ \alpha\\ \beta \epm
+\bpm
-(\alpha+\bar\alpha)\w\theta + \iC\,\eta\w\bar\eta\\
-\beta\w\theta-\alpha\w\eta\\
- \iC\,\beta\w\bar\eta - 2\iC\,\bar\beta\w\eta\\ 
-\beta\w\bar\alpha-\sigma^2_0\w\theta 
\epm,
\ee
where~$\sigma^2_0$ is now basic. From the imaginary part of~\eqref{eq: ddalpha}, this form~$\sigma^2_0$ must satisfy 
$\bigl(\overline{\sigma^2_0}\w\eta+\sigma^2_0\w\bar\eta\bigr)\w\theta=0$.  
This implies~$\sigma^2_0\w\theta = (r\eta + s\bar\eta)\w\theta$, 
where $r$~and~$s$ are real and complex functions, respectively, on~$B_3$.  
By adding~$r\theta$ to~$\sigma_0$ and calling the result~$\sigma$,
it can be arranged that $r=0$, and the structure equations become
\be
\d\bpm\theta\\ \eta\\ \alpha\\ \beta\\ \epm
=
-\bpm0&0&0&0\\ 0&0&0&0\\ \sigma&0&0&0\\ 0&\sigma&0&0\\ \epm
\bpm\theta\\ \eta\\ \alpha\\ \beta\\ \epm
+\bpm
-(\alpha+\bar\alpha)\w\theta + \iC\,\eta\w\bar\eta\\
-\beta\w\theta-\alpha\w\eta\\
- \iC\,\beta\w\bar\eta - 2\iC\,\bar\beta\w\eta
\\ 
-\beta\w\bar\alpha-s\,\bar\eta\w\theta\\ \epm
\ee
where, now, $\sigma$ is uniquely specified by these conditions.  

Thus, $B_3$ is endowed with a canonical~$\{e\}$-structure 
and this constitutes Cartan's solution of the equivalence problem. 

\subsubsection{Identities}
To complete the structure equations, 
however, a formula for~$\d\sigma$ is needed.  
The $\d(\d\alpha)=0$ equation now yields
\be
0=\left(\d\sigma-(\alpha+\bar\alpha)\w\sigma 
-\iC\,\beta\w\bar\beta\,\right)\w\theta
\ee
so that~$\d\sigma = (\alpha+\bar\alpha)\w\sigma +\iC\,\beta\w\bar\beta 
+\rho\w\theta$ where $\rho$ is a real $1$-form.  

Using this equation, the identity~$\d(\d\beta)=0$ expands to
\be
0 = \d(\d\beta) = \rho\w\eta\w\theta 
- \bigl(\d s-(3\bar\alpha+\alpha)s\bigr)\w\bar\eta\w\theta,
\ee
from which it follows that there are complex functions~$u$, $p$, and $q$ 
on~$B_3$ so that
\be
\d s = (3\bar\alpha+\alpha)s + u\,\theta + p\,\eta + q\,\bar\eta,
\ee
whence $\rho\w\theta = -(p\,\bar\eta + \bar p\,\eta)\w\theta$, so that
\be
\d\sigma = (\alpha+\bar\alpha)\w\sigma +\iC\,\beta\w\bar\beta 
-(p\,\bar\eta + \bar p\,\eta)\w\theta.
\ee

The final Bianchi identity will follow from~$\d\bigl(\d\sigma\bigr)=0$, 
and this expands to give the statement 
that there exist functions~$a$, $r$, and $v$ on~$B_3$,
with $r$ being real valued, so that
\be
\d p = (3\bar\alpha + 2\alpha)p - \iC s\,\bar\beta + a\theta + r\eta
             + v\bar\eta.
\ee

\subsection{Conclusions}\label{ssec: conclusions}
Several conclusions can be drawn from these calculations.  

\subsubsection{Automorphisms of \textrm{CR}-structures}
First of all, since the group of symmetries of a nondegenerate CR-structure 
on a~$3$-manifold embeds into the group of symmetries of an $\{e\}$-structure 
on an $8$-manifold, it follows that the group of symmetries 
of such a CR-structure is a Lie group of dimension at most~$8$.  
Moreover, this maximum dimension can be reached only if the local symmetry group 
of the $\{e\}$-structure on~$B_3$ acts with open orbits on~$B_3$.  

\subsubsection{Maximal symmetry}
However, by the structure equations, such open orbits exist 
if and only if the functions~$s$ and $p$ are locally constant.  
The structure equation for $\d s$, however, shows that~$s$ 
cannot be locally constant unless it vanishes  
(which implies, in turn, that~$p$ vanishes as well).  

In this case, the equations
\be
\begin{aligned}
\d\theta &=-(\alpha+\bar\alpha)\w\theta + \iC\,\eta\w\bar\eta\\
\d\eta &= -\beta\w\theta-\alpha\w\eta\\
\d\alpha &= -\sigma\w\theta-\iC\,\beta\w\bar\eta-2\iC\,\bar\beta\w\eta\\
\d\beta &= -\sigma\w\eta+\bar\alpha\w\beta\\
\d\sigma &= (\alpha+\bar\alpha)\w\sigma +\iC\,\beta\w\bar\beta
\end{aligned}
\ee
are the structure equations of a Lie group of dimension~$8$.  

Naturally, the reader will want to know which one.  
The simplest way to identify the group is to notice 
that there are no~$\alpha\w\bar\alpha$ terms 
on the right hand side of these equations, 
but that~$\alpha$ appears in the right hand side 
of all the equations except that of~$\d\alpha$.  
This implies that the vector 
fields~$X$ and~$Y$ dual to the real and imaginary parts 
of~$\alpha$ form a maximal torus of dimension~$2$ 
in the Lie algebra of infinitesimal symmetries of the coframing.  

For any form~$\phi$ in the coframing, define its $X$- and 
$Y$-weights by the formulae
\be
\begin{aligned}
       w_X(\phi)\,\phi &= X\lhk \d\phi,\\
\iC\,w_Y(\phi)\,\phi &= Y\lhk \d\phi.\\
\end{aligned}
\ee
Then, plotting the pairs~$\bigl(w_X(\phi), w_Y(\phi)\bigr)$ 
in the plane as~$\phi$ ranges over the 
basis~$(\theta,\eta,\bar\eta,\alpha,\bar\alpha,\beta,\bar\beta,\sigma)$ 
reveals the characteristic hexagon of the roots of~$A_2$.  
Moreover, because the roots are `half-real', 
the actual real form of~$A_2$ represented must be~$\eusu(2,1)$.  
Thus, the group must be~$\SU(2,1)$. 

In fact, when~$s$ vanishes identically, the structure equations 
can be written in the form $\d\gamma = -\gamma\w\gamma$ where
\be
\gamma = \bpm 
-{\frac13}(2\alpha+\bar\alpha)&-\iC\,\bar\beta&-\iC\,\sigma\\
\eta&{\frac13}(\alpha-\bar\alpha)&\iC\,\beta\\
-\iC\,\theta & \bar\eta & {\frac13}(\alpha+2\bar\alpha)
\epm.
\ee
Note that~$\gamma$ takes values in~$\eusu(2,1)$, where the model of~$\SU(2,1)$ 
being used is the subgroup of~$\SL(3,\bbC)$ that fixes the Hermitian form~$H$ 
in three 
variables
\be
H(Z)=Z_3\,\overline{Z_1}-Z_2\,\overline{Z_2}+Z_1\,\overline{Z_3}\,.
\ee

In particular, if~$M^3$ is simply connected, there is a smooth 
map~$F:B_3\to\SU(2,1)$ that satisfies $F^*(g^{-1}\,\d g) = \gamma$.  
As the structure equations show, $F$ maps each fiber of~$B_3\to M$ 
to a left coset of the parabolic subgroup~$P\subset\SU(2,1)$ 
consisting of the upper triangular matrices in~$\SU(2,1)$, 
i.e., the subgroup that fixes the $H$-null line~$L_0\subset\bbC^3$ 
defined by~$Z_2=Z_3=0$.  

Now, $\SU(2,1)/P$ is naturally identified with the 
hypersurface~$N^3\subset\bbC\bbP^2$ of $H$-null lines 
in~$\bbC\bbP^2$.  Thus, $F$ covers a map~$f:M^3\to N^3$ that is a local 
equivalence of CR-structures.  

The conclusion is that every CR-structure with 
$8$-dimensional infinitesimal symmetry algebra on a simply connected 
$3$-manifold  has a `developing map' to~$N^3$ 
that is unique up to composition with a CR-automorphism 
of this `flat' structure.%
\footnote{Explicitly computing this developing map requires solving a Lie 
equation of the form~$\d g = g \gamma$ where $\gamma$ is a known $1$-form 
with values in~$\eusu(2,1)$.}

\subsubsection{The non-flat case}
In the general case, $s$ is the coefficient of a tensor field that is 
well-defined on~$M$.  

The simplest such expression involving~$s$ is perhaps~$Q 
= s\bar s\,\theta^4$, 
which is a well-defined section of $S^4\bigl(D^\perp\bigr)$ 
on~$M$.  This well-definedness follows since~$Q$ is manifestly semibasic 
and a computation using the structure equations reveals that, 
for any vertical vector field~$Y$ for the projection~$B_3\to M$, 
the Lie derivative of~$Q$ with respect to~$Y$ vanishes.  

Also the expression~$S = s\,\bar\eta{\otimes}\bar\eta{\otimes}\theta$ 
can be interpreted as a well-defined section of the 
bundle~$S^{0,2}(D){\otimes}D^\perp$ over~$M$, 
i.e., the bundle of complex anti-linear quadratic forms on~$D$ 
with values in~$D^\perp$.  

Other combinations of the functions on~$B_3$ make well-defined tensors on~$M$ 
as well, but have to be treated with more care.  
For example, the expression
\be
E = s\,\bar\eta^2{\circ}\theta 
+ 2\iC p\,\bar\eta{\circ}\theta^2 \mod \theta^3
\ee 
describes a well-defined section of the quotient bundle~$S^3(T^*M)/(D^\perp)^3$. 

The verification of these statements will be left to the reader.

\subsubsection{Structure reduction in the non-flat case}
In the case where~$Q=s\bar s\,\theta^4$ is nonvanishing on~$M$, 
there is a canonical reduction of~$B_3$ to a~$\bbZ_2$-structure~$B_4\to M$ 
defined by the equations~$s=-1$, $p=0$, $u+\bar u=0$.  

This follows from the formulae for~$\d s$ and~$\d p$ 
together with the formula
\be
\d u \equiv (4\bar\alpha+2\alpha) u + p\,\beta + 4s \sigma 
\mod \theta,\eta,\bar\eta,
\ee
which is derived from the identity~$\d(\d s)=0$.  

Pulling all the given quantities back to~$B_4$,
writing~$u = 2\iC m$ where $m$ is real and 
replacing~$q$ by $8\bar q$ for notational convenience, 
this results in equations
\be
\begin{aligned}
\alpha &= \iC m\,\theta - 3q\,\eta + \bar q\,\bar\eta\\
 \beta &= \iC a\,\theta + \iC v\,\eta + \iC r\,\bar\eta,
\end{aligned}
\ee
and  in structure equations of the form
\be
\begin{aligned}
\d\theta &= 2(q\eta+\bar q\bar\eta)\w\theta +\iC\,\eta\w\bar\eta\\
\d\eta &= t\,\eta\w\theta + \bar q\,\eta\w\bar\eta -\iC r\,\bar\eta\w\theta
\end{aligned}
\ee
for some function~$t$ constructed out of the other invariants. 

Under the $\bbZ_2$-action on the double cover~$B_4\to M$, the form~$\theta$ 
is even while~$\eta$ is odd. Thus, the coframining~$(\theta,\eta)$ 
is well-defined on~$M$ up to a replacement of the form~$(\theta,\eta)\mapsto(\theta,-\eta)$.
It also follows that $t$ and $r$ are even while~$q$ is odd.

In particular, it follows from this discussion that the group of symmetries 
of a nondegenerate CR-structure for which~$Q\not=0$ 
is a Lie group of dimension at most~$3$ 
and that this upper bound is reached only for homogeneous structures, 
in which case, the functions~$q$, $r$, and $t$ must be constants.  

Indeed, if one assumes that these functions are constants, 
then computing the exterior derivatives of the above equations 
yields that $t+\bar t =0$, so that $t = \iC b$, 
for some real constant~$b$, and the equation~$rq+b\bar q = 0$.
  
Conversely, any solution~$(r,b,q)\in\bbR^2\times\bbC$ of~$rq+b\bar q = 0$ 
defines a homogeneous CR-structure.  Cartan used this fact in his 1932 
papers to classify the homogeneous CR-hypersurfaces in~$\bbC^2$.

\section{Generalizations}\label{sec: generalizations}

In this brief section, some remarks will be made about 
generalizations of these constructions to higher dimensions.

The reader may wonder just how much of this geometry
for real hypersurfaces in unimodular complex surfaces 
can be generalized to higher dimensions.  In fact, 
a great deal of it can be.  Such generalizations are
not the focus of this article, but the following comments
may be a useful guide to the reader who wants to explore
this further.

As the reader will no doubt realize, the geometric 
structures induced on hypersurfaces discussed here 
are quite different from the ones considered by
Huisken and Klingenberg~\cite{MR01f:53141,MR03d:32046}.
Consequently, the corresponding flows are unrelated
to the flows that they consider.

\subsection{Hypersurfaces in unimodular complex manifolds}
\label{ssec: hypsurfin1modmfds}
Let~$n\ge 1$ be fixed
and let~$X$ be a complex manifold of dimension~$n{+}1$
and~$\Upsilon$ be a nowhere-vanishing holmorphic~$(n{+}1,0)$-form
on~$X$.  The pair~$(X,\Upsilon)$ will be referred to as a
\emph{unimodular complex $(n{+}1)$-manifold}.

Let~$M^{2n+1}\subset X$ be a real hypersurface that is smoothly
embedded in~$X$.  (The immersed case does not differ substantially
and lower regularity will suffice for the constructions to be 
carried out in this subsection, but these refinements will be 
left to the interested reader.)  Let~$\Phi = M^*\Upsilon$ be
the pullback of~$\Upsilon$ to~$M$.  

By linear algebra, there exists, at least locally on~$M$, a
nonvanishing $1$-form~$\theta=\overline{\theta}$ such 
that~$\theta\w\Phi=0$. This $1$-form is unique up to a real
multiple. The kernel of~$\theta$ is a codimension~$1$ plane
field~$D\subset TM$ that is, by definition the complex tangent
space of~$M$.  One says that~$M$ is \emph{nondegenerate} if $\d\theta$
is a nondegenerate $2$-form on~$D$, i.e., if~$\theta\w(\d\theta)^n$
is nowhere-vanishing on the domain of definition of~$\theta$.
In any case, $\d\theta$ restricts to~$D$ to be a real $(1,1)$-form.
(Of course, this is, up to multiples, the Levi form of~$M$.)
One says that~$M$ is \emph{strictly pseudoconvex} if it is
possible to choose~$\theta$ such that~$\d\theta$ is a positive
definite~$(1,1)$-form on~$D$.  This latter requirement, if 
satisfiable, determines~$\theta$ up to a positive multiple.

For simplicity, only the strictly pseudoconvex case will be
discussed any further here.  The generalization of 
Proposition~\ref{prop: cancoframing} to all cases~$n\ge1$ is then

\begin{proposition}[The canonical $\SU(n)$-structure]
\label{prop: cancoframinggenn}
Let~$M\subset X$ be a strictly pseudoconvex real hypersurface
in a unimodular complex manifold~$(X^{n+1},\Upsilon)$ and let~$\Phi
= M^*\Upsilon$. Then there is a unique $\SU(n)$-structure on~$M$ 
whose sections~$(\theta,\eta_1,\ldots,\eta_n)$ 
satisfy~$\theta=\overline\theta$ and
\begin{enumerate}
\item $\Phi = \theta\w\eta_1\w\ldots\w\eta_n$\,,
\item $\d\theta = \iC\bigl(\eta_1\w\overline{\eta_1}
                       +\cdots+\eta_n\w\overline{\eta_n}\bigr)$\,.
\end{enumerate}
\end{proposition}

\begin{proof}
The proof follows exactly the same lines as that 
of Proposition~\ref{prop: cancoframing}, 
so there is no need to repeat it here.
\end{proof}

Of course, Proposition~\ref{prop: cancoframinggenn} implies that there
is a canonical normal associated to a strictly pseudoconvex real
hypersurface, so there will be an evolution equation for such hypersurfaces
generalizing the $n=1$ case.  

There is also a canonical volume form
defined as~$\frac{2^{-n}}{n!}\theta\w(\d\theta)^n$, 
which is the volume form associated to the canonical metric
\be
\d s^2 = \theta^2 
+ \eta_1\circ\overline{\eta_1} + \cdots + \eta_n\circ\overline{\eta_n}.
\ee

As in the $n=1$ case, the primary invariants are found by computing the
derivatives of the remaining forms in a coframing of the $\SU(n)$-structure.
Computation (i.e., expanding the identities~$\d(\d\theta)=
\d(\theta\w\eta_1\w\cdots\w\eta_n)=0$
and applying exterior algebra identities) yields that, 
for such a coframing~$(\theta,\eta_1,\ldots,\eta_n)$, 
there are structure equations
\be\label{eq: genn1ststreqs}
\begin{aligned}
\d\theta &= \iC\bigl(\eta_1\w\overline{\eta_1}
                       +\cdots+\eta_n\w\overline{\eta_n}\bigr)\,,\\
\d\eta_j &= -\phi_{j\bar k}\w\eta_k 
             +2\iC\,\theta\w\bigl(\,a\,\eta_j+b_{jk}\,\overline{\eta_k}\bigr),
\end{aligned}
\ee
for some unique real-valued function~$a = \bar a$ (independent 
of the choice of coframing), some unique complex functions~$b_{ij}=b_{ji}$, 
and some unique $1$-forms~$\phi_{j\bar k}=-\overline{\phi_{k\bj}}$ 
satisfying $\phi_{1\bar1}+\cdots+\phi_{n\bar n}=0$. 

The invariant $a$, which is fourth order, 
has the same `mean curvature' interpretation relative to the canonical 
volume form as it does in the~$n=1$ case. 
 
The quadratic expression~$B = \overline{b_{jk}}\,\eta_j\circ\eta_k$
is also fourth order and well-defined, independent of choice of 
coframing.  When~$B$ vanishes identically, one finds that~$a$ is
constant.  Moreover, the quadratic form~$\d\sigma^2 
= \eta_1\circ\overline{\eta_1}+\cdots+\eta_n\circ\overline{\eta_n}$
descends to the leaf space~$Z$ of the Reeb vector field~$T$ on~$M$ to define
a K\"ahler-Einstein metric on~$Z$ with associated K\"ahler form~$\frac12\d\theta$
and Ricci form~$\Ric(\d\sigma^2) = 4na\,\d\sigma^2$.  
When~$a\not=0$, the original~$M$ can be recovered (up to a covering)
as the circle bundle of $(n,0)$-forms of a fixed norm with
respect to~$\d\sigma^2$, regarded as a real hypersurface in the
total space of the canonical bundle of~$Z$ (endowed with its tautological
holomorphic $(n{+}1)$-form).  When~$a=0$ and~$h$ is a (local) K\"ahler
potential for the Ricci-flat K\"ahler metric~$\d\sigma^2$, 
one can recover~$M$ locally as the hypersurface in~$\bbC\times Z$ defined 
by~$\text{Im}(z^0) = h$, where the holomorphic~$(n{+}1)$-form
on~$\bbC\times Z$ is~$\Upsilon=\d z^0\w\Psi$ and~$\Psi$ 
is a (locally defined) $\d\sigma^2$-parallel 
holomorphic volume form on~$Z$.

Just as in the case~$n=1$, the function~$a$
and the tensor~$B$ are all of the fourth order invariants 
of the strictly pseudoconvex hypersurface under the action 
of the unimodular biholomorphism pseudogroup.  However, when~$n>1$,
these invariants and the invariants derived from them 
do not constitute a complete set of invariants.

The 1-forms~$\phi_{j\bar k}$ are the connection forms (relative
to the chosen coframing) of a canonical connection on the~$\SU(n)$-structure.
To generate a complete set of invariants in Cartan's sense,
one must include, along with~$a$ and~$B$ and their covariant
derivatives, the curvature tensor of the connection~$\phi = (\phi_{j\bar k})$ 
(a tensor that is of fifth order) and its covariant derivatives.

The $\SU(n)$-structure and its canonical connection
constitute the solution of the equivalence problem
(in Cartan's sense) for strictly pseudoconvex real hypersurfaces 
in unimodular complex manifolds. Higher order invariants can be defined 
and studied by differentiating these equations.

There is also a generalization of the embedding results derived
in the~$n=1$ case:  Every real-analytic $\SU(n)$-structure on~$M^{2n+1}$ 
that satisfies the first order structure equations~\eqref{eq: genn1ststreqs}
(i.e., whose intrinsic torsion has this form) arises as the
canonical $\SU(n)$-structure associated to a real-analytic,
strictly pseudoconvex real hypersurface in a unimodular 
complex manifold of dimension~$n{+}1$ and this realization is
essentially unique.  The proof is a straightforward
application of the Cartan-K\"ahler Theorem and will be left 
to the reader.

\subsection{Hypersurfaces in complex manifolds with
a specified volume form}
\label{ssec: specificvolform}
Finally, it should be pointed out that one can even generalize
the construction to cover the case real hypersurfaces in 
a complex $(n{+}1)$-manifold~$X$ endowed with a positive, 
real $(n{+}1,n{+}1)$-form, i.e., a volume form~$\Omega$ 
(in the usual sense) on~$X$.  

The point is that, for a real hypersurface~$M\subset X$ 
with definite Levi form, one can use~$\Omega$ to 
make a canonical choice of~$\theta$ as a real-valued 
$1$-form on~$M$ whose kernel is the complex tangent space
of~$M$ and whose associated volume form~$\theta\w(\d\theta)^n$
is determined by~$\Omega$.  

More precisely, let~$\theta_0$ be any $1$-form on an open 
set~$U\subset M$ whose kernel is the complex tangent 
bundle~$D\subset TM$ and such that~$\d\theta_0$ restricts to~$D$
to be a positive $(1,1)$-form.  These conditions determine~$\theta_0$
up to multiplication by a positive function~$f$ on~$U$.  Let~$T_0$
be the Reeb vector field on~$M$ associated to~$\theta_0$, i.e.,
$\theta_0(T_0)=1$ and~$\d\theta_0(T_0,Y)=0$ for all tangent vector
fields~$Y$ on~$U$.  Let~$f$ be any positive smooth function on~$U$.
Then the Reeb vector field~$T$ associated to~$\theta = f\,\theta_0$
is of the form~$T = (1/f)\,T_0 + S$ where~$S$ is a section of~$D$.
Since~$\theta\w(\d\theta)^n = f^{n+1}\,\theta_0\w(\d\theta_0)^n$
and since~$M^*\bigl((\iC T)\lhk\Omega\bigr) 
= f^{-1}\,M^*\bigl((\iC T_0)\lhk\Omega\bigr)$, it follows that there
is a unique choice of~$f>0$ such that~$\theta\w(\d\theta)^n
= M^*\bigl((\iC T)\lhk\Omega\bigr)$.  Once~$\theta$ is fixed, 
its Reeb vector field~$T$ is determined, yielding a canonical
normal vector field~$\iC T$ along~$M$, which can then be used to define 
an evolution for such hypersurfaces that is invariant under the
pseudogroup of biholomorphisms that preserve~$\Omega$.

However, the reader should be aware that this is, in some
sense, not much of a generalization.  In the first place,
the pseudogroup of local biholomorphisms of~$X$ that preserve~$\Omega$
will, generally, be `smaller' than the unimodular
biholomorphism group because the pseudogroup of biholomorphisms 
that preserves~$\Omega$ must also preserve 
the~$(1,1)$-form~$\Ric(\Omega)$.  Usually, this~$(1,1)$-form will
be nonzero, forcing the pseudogroup to preserve a
(usually singular) complex foliation.  On the other hand, 
when~$\Ric(\Omega)$ vanishes identically, the pseudogroup perserves
a holomorphic volume form up to a constant multiple, so the
pseudogroup is only slightly larger than the unimodular biholomorphism
pseudogroup and the picture of the invariants is essentially
the same as in the unimodular case.  In the second place, 
the analysis of the local invariants is rendered somewhat more
complicated because the geometry of the form~$\Ric(\Omega)$
must be taken into account in writing down local normal forms
and computing evolution identities.  Most likely, this complication
is more of a nuisance than a serious difficulty in studying
the geometry of the evolution equation in this more general 
case, but this will be left for the interested reader to pursue.

\bibliographystyle{hamsplain}

\providecommand{\bysame}{\leavevmode\hbox to3em{\hrulefill}\thinspace}

\end{document}